# Oscar Sheynin

# Theory of Probability

# An Elementary Treatise against a Historical Background

## Contents

# 0. Introduction
## 0.1. Some Explanation

This treatise is written on an elementary level; in more difficult cases the final formulas are provided without proof. Nevertheless, it was impossible to leave out integrals and I also had to differentiate an integral with respect to a parameter. I include many examples taken from the history of probability and hope that my subject has thus become lively. I especially quote Karl Pearson's (1978, p. 1) repentant confession:

*I do feel how wrongful it was to work for so many years at statistics and neglect its history.*

In spite of a few mistakes, his book deserves serious attention.
Thus, in § 4.1.1 I criticize his opinion about Jakob Bernoulli.

I have devoted much attention to the notion of probability which fully conforms to Langevin's statement (1913/1914, p. 3):

*Dans toutes ces questions* [in the kinetic theory] *la difficulté principale est, comme nous le verrons, de donner une définition correcte et claire de la probabilité.*

Note however that *correct definition* sounds strangely.

## 0.2. The Object of the Theory of Probability

Toss a coin and the outcome will be either heads or tails. Toss it 50 times and theoretically there will be 0, 1, 2, …, 49 or 50 heads. For the time being I emphasize that the number of heads will only be determined stochastically (= probabilistically): there will be from 20 to 30 heads with such-and-such probability; from 22 to 28 heads with another probability etc. In probability theory, it will never be possible to provide a quite definite answer whereas, for example, the number of roots of a given algebraic equation can be stated at once.

Games of chance (of which coin tossing is an example) was the main subject of the early theory of probability. Their outcome depends on chance rather than on the gamblers' skill, and even now they are methodically (and therefore pedagogically as well) interesting.

Many tosses provide an example of mass random events which occur in most various settings: in population statistics (births, marriages, deaths), when treating numerous observations corrupted by unavoidable random errors, applying acceptance sampling of manufactured articles with a stochastic estimation of its error, and in various branches of knowledge (kinetic theory, epidemiology etc). And so, *the theory of probability studies mass random events, or, more precisely, their regularities* which really exist. An isolated event is random, but a (homogeneous) set of events displays regularities. Aristotle (*Metaphysics* 1026b) remarked that *none of the traditional*

*sciences busies itself about the accidental*. As stated above, neither does the theory of probability!

Laplace quite successfully applied probability to studying mass random events, and thus to investigating laws of nature (especially astronomy) and population statistics. And unlike his predecessors, he regarded the theory of probability as a branch of applied mathematics (and separated himself from the geometers: *let the geometers study …*).

I ought to add the reasonable but indefinite Laplace's opinion (1814/1886, p. CLIII): *La théorie des probabilités n'est, au fond, que le bon sens réduit au calcul*. He did not mention mass random events and furthermore his definition pertained to mathematics of his time as a whole.

Times change and we change with time … A mathematically formulated definition of the aims of the theory of probability became needed, and Boole (1851/1952, p. 251) provided it, in a seemingly dull wording: *Given the separate probabilities of any* [logical] *proposition, to find the probability of another proposition*. A similar statement pertaining to events was due to Chebyshev (1845/1951, p. 29): the theory of probability *has as its subject the determination of an event given its connection with events whose probabilities are given*. He added that probability signifies *some magnitude subject to measurement*. Prokhorov & Sevastianov (1999, p. 77) confirmed that aim and noted that such determinations were possible owing to the stability of those same mass random *phenomena*, as they stated. Anyway, owing to the stability of statistical probability (§ 1.1.3).

Since the theory of probability is axiomatized, it belongs to pure mathematics rather than a branch of applied mathematics (Laplace, see above).

## Chapter 1. Main Notions, Theorems and Formulas
### 1.1. Probability

**1.1.1.** *Theoretical Probability.* Suppose that the outcome of a trial depends on $n$ incompatible and equally possible cases only $m$ of which are favourable for the appearance of some event $A$. Then its probability is assumed as

$$P(A) = m/n \qquad (1.1)$$

and it can change from 0 to 1, from an impossible to a certain event. This is the so-called classical definition due (not to Laplace, but) to De Moivre (1711/1984, p. 237) although he formulated it in the language of chances as he also did later, in 1756.

That definition had been known or intuitively applied from antiquity. The Talmud recognized seven *levels* of food containing differing relative amounts of a prohibited element (Rabinovitch 1973, p. 41). In the $14^{th}$ century, Oresme (1966, p. 247) possibly thought about probability in the modern way since he stated that *two* [randomly chosen] *numbers were probably incommensurable*. (For us, his understanding of incommensurability was unusual.) The same idea of probability is seen in Kepler (1596). Finally, I cite Jakob Bernoulli

(1713, Chapter 1 in Pt. 4). He introduced probability just as De Moivre did but not formally, nor did he apply it in the sequel.

In ancient times, geometry started by definitions of a point, a line and a plane. The point, for example, was something dimensionless. Nowadays, such *negative* definitions are unacceptable; just consider: a man is not a woman … and a woman is not a man! We have to accept such initial notions without defining them. Here is Markov (1900; 1908, p. 2; 1924, p. 2; and 1911/1981, p. 149):

*Various concepts are defined not by words, each of which can in turn demand definition, but rather by* [our] *attitude towards them ascertained little by little.*

*I shall not defend these basic theorems linked to the basic notions of the calculus of probability, notions of equal probability, of independence of events, and so on, since I know that one can argue endlessly about the basic principles even of such a precise science as geometry.*

Then, Kamke (1933, p. 14) noted: *Um das Jahr 1910 konnte man in Göttingen das Bonmot hören*:

*Die mathematische Wahrscheinlichkeit ist ein Zahl, die zwischen Null und Eins liegt und über die man sonst nicht weis.*

At that tine, Göttingen was considered the mathematical world centre, but in 1934 Hilbert, who had been working there, stated that after the Jewish scholars were ousted, the university ceased to exist. Not without reason Khinchin (1961/2004, p. 396) noted that

*Each author* […] *invariably reasoned about equally possible and favourable chances, attempting, however, to leave this unpleasant subject as soon as possible.*

Indeed, definition (1.1) begs the question: probability depends on equal possibilities, that is, on equal probabilities. More important, it is not really a definition, but only a formula for calculating it. Just the same, the area of a square can be calculated, but the appropriate formula does not tell us the meaning of *area*. And, finally, equal possibilities exist rarely so that the application of formula (1.1) is severely restricted.

In accord with Hilbert's recommendation (1901, Problem No. 6), the contemporary theory of probability is axiomatic, but in practice statistical probability (see § 1.1.3) reigns supreme.

*Example* (application of theoretical probability). Apparently during 1613 – 1623 Galileo wrote a note about a game with 3 dice first published in 1718 (David 1962, pp. 65 – 66; English translation, pp. 192 – 195). He calculated the number of all the possible outcomes (therefore, indirectly, the appropriate probabilities) and compared the appearance of 9 or 12 points and 10 or 11 points (events *A* and *B*). Both *A* and *B* occurred in six ways; thus, *A* can appear when the number of points on the dice is (3, 3, 3) or (1, 4, 4 or 2, 2, 5) or (1, 2, 6 or 1, 3, 5 or 2, 3, 4), i. e. when the number of points on each die is the same; when it is only the same on two dice; and when it is different. However, the first case is realized only once, the second case, in 3 ways, and the last one, in 6 ways. Event *A* therefore appears in 25 ways whereas event *B*, according to similar considerations, in 27 ways.

The total number of possible outcomes is 216, 108 for 3, 4, …, or 10 points, and again 108 for 11, 12, …, 18 points and the probabilities of *A* and *B* are 25/216 and 27/216.

This example is instructive: it shows that the cases in formula (1.1) if unequally likely can be subdivided into equally possible ones. Galileo also stated that gamblers knew that *B* was more advantageous than *A*. They could have empirically compared not 25/216 and 27/216, but 25/52 and 27/52 by only paying attention to the two studied events.

*Some definitions*. When two events, *A* and *B*, have occurred, we say that their product *AB* had appeared. When at least one of them has occurred, it was the appearance of their sum, (*A* + *B*), and if only one (say, *A* but not *B*), then it was their difference (*A* – *B*).

*Example*. Two chess tournaments are to be held. The probabilities of a certain future participant to win the first place (events *A* and *B*) are somehow known. If the tournaments will occur at the same time, the product *AB* is senseless, formula (1.1) cannot be applied, the probability of that product does not exist.

**1.1.1.-1.** *The addition theorem*. For incompatible events *A* and *B*

$$P(A + B) = P(A) + P(B).$$

*Examples*. Suppose that an urn contains *n* balls, *a* of them red, *b*, blue, and *c*, white. Required is the probability of drawing a coloured ball (Rumshitsky 1966). The answer is obviously *a/n* + *b/n*.

Here, however, is only a seemingly similar problem. A die is rolled twice. Required is the probability that one six will appear. Call the occurrence of 6 points in the first and the second trial by *A* and *B*. Then

$$P(A) = 1/6, P(B) = 1/6, P(A + B) = 1/3.$$

But something is wrong! After 6 trials the probability will be unity, and in 7 trials?.. The point is, that *A* and *B* are not incompatible. See the correct solution in § 1.1.1-2.

The addition and the multiplication (see below) theorems for intuitively understood probabilities have actually been applied even in antiquity. Aristotle (*De Caelo* 292a30 and 289b22) stated that

*Ten thousand Coan throws* [whatever that meant] *in succession with the dice are impossible and it is therefore difficult to conceive that the pace of each star should be exactly proportioned to the size of its circle*.

Imagined games of chance had illustrated impossible events: the stars do not rotate around the sky randomly. Note that the naked eye sees about six thousand stars.

**1.1.1-2.** *Generalization: the formula of inclusion and exclusion*. For two events *A* and *B* the general addition formula is

$$P(A + B) = P(A) + P(B) - P(AB).$$

Indeed, in the example in § 1.1.1-1 *P*(*AB*) = 1/36 and

$$P(A + B) = 1/6 + 1/6 - 1/36 = 11/36.$$

The number of favourable cases was there 11 rather than 12. For 3 events we have

$$P(A + B + C) = P(A) + P(B) + P(C) - P(AB) - P(AC) - P(BC) + P(ABC)$$

and in the general case

$$P(A_1 A_2 \ldots A_n) = P(\sum_i A_i) - P(\sum_{i<j} A_i A_j) + P(\sum_{i<j<k} A_i A_j A_k) - \ldots$$

This *formula of inclusion and exclusion* was applied by Montmort (1708). It is a particular case of the proposition about the mutual arrangement of arbitrary sets. The conditions $i < j$, $i < j < k$, … ensure the inclusion of all subscripts without repetition. Thus, for 4 events $i < j$ means that allowed are events with subscripts 1, 2; 1, 3; 1, 4; 2, 3; 2, 4 and 3, 4, six combinations in all (of 4 elements taken 2 at a time).

**1.1.1-3.** *The multiplication theorem.* We introduce notation $P(B/A)$, denoting the probability of event $B$ given that event $A$ had occurred. Now, the theorem:

$$P(AB) = P(A)P(B/A). \qquad (1.2)$$

Switch $A$ with $B$, then

$$P(AB) = P(B)P(A/B). \qquad (1.3)$$

*Example* 1 (Rumshitsky 1966). There are 4% defective articles in a batch; among the others 75% are of the best quality. Required is the probability that a randomly chosen article will be of the best quality.

Denote the extraction of a standard article by $A$, and by $B$, of one of the best. Then

$$P(A) = 1 - 0.04 = 0.96;\ P(B/A) = 0.75.\ P(AB) = 0.96 \cdot 0.75 = 0.72.$$

*Example* 2. What number of points, 11 or 12, will occur more probably in a cast of two dice? Leibniz (Todhunter 1865, p. 48) thought that both outcomes were equally probable since each was realized in only one way, when casting 5 and 6 and 6 and 6 respectively. An elementary mistake committed by a great man! Denote by $A$ and $B$ the occurrence of 5 and 6 on a die, then $P(A) = 1/6$, $P(B) = 1/6$. Yes, both outcomes after casting both dice are the same, $P(AB) = P(A)P(B) = 1/36$, but we ought to take into account that the first alternative can appear in two ways, 5 and 6, and 6 and 5, and is therefore twice as probable.

In general, if $P(B/A) = P(B)$ the multiplication theorem is written as

$$P(AB) = P(A)P(B),$$

and the events $A$ and $B$ are called *independent*.

*Example* 3. $A$ and $B$ have 12 counters each and play with 3 dice. When 11 points appear, $A$ gives $B$ a counter and $B$ gives a counter to $A$ when 14 points occur. Required are the gamblers' chances of winning all the counters. This is Additional problem No. 5 formulated by Pascal, then by Huygens (1657), who provided the answer without solution.

There are 216 outcomes of a cast of 3 dies, 15 of them favouring the appearance of 14 points, and 27 favouring 11 points, see Example in § 1.1.1. The probabilities or chances of winning are therefore as 15/27 = 5/9. For winning 12 counters the chances therefore are as $5^{12}/9^{12}$.

This was the first of a series of problems describing the *gambler's ruin*. They proved extremely interesting and among their investigators were De Moivre and Laplace. In a particular case, the fortune of one of the gamblers was supposed to be infinite.

A series of games of chance can be thought of as a random walk whereas, when considered in a generalized sense, they become a random process (§ 5.2).

Suppose now that more than 2 (for example, 4) events are studied. The multiplication theorem will then be generalized:

$$P(A_1A_2A_3A_4) = P(A_1)P(A_2/A_1)P(A_3/A_1A_2)P(A_4/A_1A_2A_3).$$

The last multiplier, for example, denotes the probability of event $A_4$ given that all the other events had happened.

Reader! Bear with me for some time yet; two more statements are needed, perhaps not very elegant (every man to his taste).

**1.1.1-4.** *A more essential generalization of the multiplication theorem.* Suppose that event $A$ can occur with one and only one of several incompatible events $B_1, B_2, …, B_n$. It follows that our notation $P(AB)$ can now be replaced simply by $P(A)$, so that formula (1.3) will be

$$P(A) = P(B_1) P(A/B_1) + P(B_2) P(A/B_2) + … + P(B_n) P(A/B_n) =$$

$$\sum_{i=1}^{n} P(B_i)P(A/B_i). \qquad (1.4)$$

This is the formula of *total probability* and the $B_i$'s may be considered the *causes* of the occurrence of $A$, each leading to $A$ although only with its own probability.

Suppose that 3 urns have, respectively, 1 white (w) and 2 black (b) balls; 2 w and 1 b ball; and 3 w and 5 b balls. An urn is randomly selected and a ball is drawn from it. Required is the probability that that ball is white.

The probabilities of extracting a white ball from those urns are $P(A/B_i) = 1/3, 2/3$ и $3/8$, and the probability of selecting any urn is the same, $P(B_i) = 1/3$. Therefore,

$$P(A) = (1/3)·(1/3) + (1/3)·(2/3) + (1/3)·(3/8) = 0.458 < 1.$$

It is quite possible that the extracted ball was black. But can $P(A) > 1$? Let

$$P(A/B_1) = P(A/B_2) = P(A/B_3) = 0.99$$

and of course

$$P(B_1) + P(B_2) + P(B_3) = 1.$$

But the feared event will not occur even when the first 3 probabilities are so high. But can $P(A) = 1$?

**1.1.1-5.** *The Bayes formula.* The left sides of equations (1.2) and (1.3) coincide, and their right sides are equal to each other

$$P(A)P(B/A) = P(B)P(A/B),$$

or, in previous notation,

$$P(A)P(B_i/A) = P(B_i)P(A/B_i),$$

$$P(B_i/A) = \frac{P(B_i)P(A/B_i)}{P(A)}. \qquad (1.5)$$

Replace finally $P(A)$ according to formula (1.4):

$$P(B_i/A) = \frac{P(B_i)P(A/B_i)}{\sum_{i=1}^{n} P(B_i)P(A/B_i)}. \qquad (1.6)$$

It is time to contemplate. We assigned probabilities $P(B_1)$, $P(B_2)$, …, $P(B_n)$ to causes $B_1$, $B_2$, …, $B_n$ and they are in the right side of (1.6). But they are prior whereas the trial was made: the event $A$ has occurred and those prior probabilities can now be corrected, replaced by posterior probabilities $P(B_1/A)$, $P(B_2/A)$, …, $P(B_n/A)$.

Bayes (1764) included formula (1.6) but only in the particular case of $n = 1$ (which means going back to the previous formula). However, it is traditionally called after him. More precisely, from 1830 onwards it was formula (4.5) that was called after him. Nevertheless, Cournot (1843, § 88), although hesitantly, attributed formula (1.6) to Bayes; actually, it appeared in Laplace's great treatise (1812, § 26).

And who was Bayes? A talented mathematician. His posthumous memoir (1764 – 1765) became lively discussed in the early 20[th] century since prior probabilities were rarely known; is it possible to suppose that they are equal to each other? Laplace (1814/1995, p. 116) thought that hypotheses should be created without attributing them *any reality* and *continually* corrected by new observations. Discussions are still continuing and anyway many terms are called after Bayes, for example, Bayesian approach, estimator etc.

*Example.* Consider the same three urns as above. For them, the fractions in the right side of formula (1.5) differ one from another only by multipliers $P(A/B_i)$, which are to each other as $(1/3):(2/3):(3/8) = 8:16:9$. The same can therefore be stated about the posterior probabilities $P(B_i/A)$. It is certainly possible to take into consideration the previously established value $P(A) = 0.458$ and calculate them:

$$P(B_1/A) = \frac{1}{3 \cdot 3 \cdot 0.458}, \quad P(B_2/A) = \frac{2}{3 \cdot 3 \cdot 0.458}, \quad P(B_3/A) = \frac{3}{3 \cdot 8 \cdot 0.458}.$$

Understandably, probability $P(B_2/A)$ turned out as the highest of them: the relative number of white balls was largest in that same urn, the second one.

Stigler (1983/1999) applied the *Bayes theorem* in the mentioned particular case for stating that another English mathematician, Saunderson, was the real author of the Bayes memoir. He (p. 300) assigned subjective probabilities to three differing assumptions (for example, did each of them, Bayes and Saunderson, keep in touch with De Moivre) and multiplied these probabilities for each of the two. Their ratio occurred to be 3:1 in favour of the latter. Tacitly allowing an equality of the corresponding prior probabilities, Stigler (p. 301) decided that the probability of Saunderson's authorship was three times higher. Stigler's tacit assumption was absolutely inadmissible and that his (happily forgotten) conclusion ought to be resolutely rejected. That same Stigler allowed himself to vomit an abuse on Euler (§ 6.2) and Gauss (Sheynin 1999a, pp. 463 – 466).

**1.1.1-6.** *Subjective probability*. It is naturally prior and somewhat complements the theoretical probability (1.1). Indirectly, it is applied very often, especially when there exists *no reason* to doubt the existence of equal probabilities of some outcomes. Thus, the probability of each outcome of a cast of a die is supposed to be 1/6, although any given die is not exactly *regular*. Poisson and Cournot (1843/1984, p. 6) were the first to mention it definitely. They even called it and the objective probability by different terms, *chance* and *probability*.

Here is Poisson's (1837, § 11) instructive problem. An urn contains $n$ white and black balls in an unknown proportion. Required is the probability that an extracted ball is white. The number of white balls can be 0, 1, 2, ..., $n$, – ($n$ + 1) allegedly equally probable cases. The probability sought is therefore the mean of all possible probabilities

$$\frac{1}{n+1}(\frac{n}{n}+\frac{n-1}{n}+...+\frac{1}{n}+\frac{0}{n}) = \frac{1}{2}$$

*as it should have been*. His answer conforms to the principles of the information theory which Poisson himself understood perfectly well: it, his answer, corresponded to *la perfaite perplexité de notre esprit*.

Poisson (1825 – 1826) applied subjective probability when investigating a game of chance. Cards are extracted one by one from six decks shuffled together as a single whole until the sum of the points in the sample obtained was in the interval [31; 40]. The sample is not returned and a second sample of the same kind is made. It is required to determine the probability that the sums of the points are equal. Like the gamblers and bankers, Poisson tacitly assumed that the second sample was extracted as though from the six initial fresh decks.

Actually, this was wrong, but the gamblers thought that, since they did not know what happened to the initial decks, the probability of drawing some number of points did not change.

When blackjack is played, bankers are duty bound to act the same wrong way: after each round the game continues without the used cards, and, to be on the safe side, they ought to stop at 17 points. A gambler endowed with a retentive memory can certainly profit from this restriction.

Here are other examples. *Redemption of the first born*. The Jerusalem Talmud (Sanhedrin 1[4]) describes how lots were taken. The main point was that the voters were afraid that there will be no more special ballots left freeing the last voters from the payment. They actually thought about the subjective probabilities of the distribution of those special ballots among consecutive voters. Tutubalin (1972, p. 12) considered the same problem in quite another setting and proved that the fears of the voters were unfounded.

*Another example.* Rabinovitch (1973, p. 40) desribed the statement of Rabbi Shlomo ben Adret (1235 – 1310 or 1319) about eating some pieces of meat one of which was not kosher. One piece after another may be eaten because (actually) the probability of choosing the forbidden piece was low, and when only two pieces are left, – why, the forbidden piece was already eaten, so eat these two as well!

Subjective *opinions* are mathematically studied, for example those pertaining to expert estimates and systems of voting. In those cases the merits of the economic projects or candidates are arranged in ascending or descending order of preference, see § 5.1.

**1.1.2.** *Geometrical Probability.* The classical definition of probability can be generalized, and, in a manuscript of 1664 – 1666, Newton (1967, pp. 58 – 61) was the first to do so. He considered a ball falling upon the centre of a circle divided into sectors whose areas were in *such proportion as* 2 *to* √5. If the ball *tumbles* into the first sector, a person gets *a*, otherwise he receives *b*, and his *hopes is worth*

$$(2a + b\sqrt{5}) \div (2 + \sqrt{5}).$$

The probabilities of the ball tumbling into these sectors were as 2 to √5, as Newton also indirectly stated. See also Sheynin (1971a).

The classical definition is still with us with *m* and *n* being real rather than only natural numbers. In this way many authors effectively applied *geometrical probability*. Buffon (1777, § 23) definitively introduced it by solving his celebrated problem. A needle of length 2*r* falls *randomly* on a set of parallel lines. Determine the probability *P* that it intersects one of them. It is easily seen that

$$P = 4r/\pi a$$

where *a* > 2*r* is the distance between adjacent lines. Buffon himself had however only determined the ratio *r/a* for *P* = 1/2. His main aim was to *mettre donc la Géométrie en possession de ses droits sur la science du hazard* (Buffon 1777/1954, p. 471). Later authors generalized the Buffon problem, for example by replacing lines by rectangles or squares.

Laplace (1812, chapter 5) noted that after, say, 100 such trials the number π can be calculated. He thus suggested the Monte Carlo method (of statistical simulation). A formal definition of the new concept was only due to Cournot (1843, § 18). More precisely, he offered a general definition for a discrete and a continuous random variable by stating that probability was the ratio of the *étendue* of the favourable cases to that of all the cases. We would now replace *étendue* by *measure* (in particular, by area).

Actually, beginning with Nikolaus Bernoulli (1709/1975, pp. 296 – 297), see also Todhunter (1865, pp. 195 – 196), each author dealing with continuous laws of distribution (§ 2.1) applied geometric probability. The same can be said about Boltzmann (1868/1909, p. 49) who defined the probability of a system being in a certain phase as the ratio of the time during which it is in that time to the whole time of the motion. Ergodic theorems can be mentioned, but they are beyond our boundaries.

Determine the probability that a random chord of a given circle is shorter than the side of an inscribed equilateral triangle (Bertrand 1888). This celebrated problem had been discussed for more than a century and several versions of *randomness* were studied. Bertrand himself offered three different solutions, and it was finally found out that, first, there was an uncountable number of solutions, and, second, that the proper solution was *probability equals* 1/2 which corresponded to *la perfaite perplexité de notre esprit* (§ 1.1.1-6).

Finally, the encounter problem (Laurent 1873, pp. 67 – 69): two persons are to meet at a definite spot during a specified time interval (say, an hour). Their arrivals are independent and occur *at random*; the first one to come waits only for a certain time (say, 20 minutes), then leaves. Required is the probability of a successive encounter.

Denote the time of their arrivals by *x* and *y*, then $|x – y| \leq 20$ or $|y – x| \leq 20$, and a graphical solution is simple and instructive, see also § 3.2.

**1.1.3.** *Statistical Probability.* Suppose that a random event occurred μ times in ν trials. Then its relative frequency (frequency, as I will call it) or statistical probability is

$$\hat{p} = \mu/\nu \qquad (1.7)$$

and it obviously changes from 0 to 1.

Newton (§ 1.1.2), while commenting on his second thought experiment, a roll of an irregular die, concluded that, nevertheless, *It may be found how much one cast is more easily gotten then another*. He likely had in mind statistical probabilities rather than analytic calculations. And he may well have seen Graunt's pioneer statistical contribution of 1662 where all deductions pertaining to population and medical statistics had been based on statistical probabilities.

Statistical probability was applied even by Celsus (1935, p. 19) in the first century of our era:

*Careful men noted what generally answered the better, and then began to prescribe the same for their patients. Thus sprang up the Art of medicine.*

He certainly had no numerical data at his disposal, but qualitative statements had been a distinctive feature of ancient science.

The definition above is only meaningful if the trials are mutually independent and the calculated probability remains almost the same in a subsequent series of similar trials. If results of some trials essentially differ, say, from one day of the week to another, then each such day ought to be studied separately. And what kind of trials do we call independent? For the time being, we say: trials, whose results do not influence each other, also see § 1.1.1-3.

The imperfection of the theoretical probability and its narrow field of applications led to the appearance of the statistical probability as the main initial notion (Richard Mises, in the 1920s).

A rigorous implementation of his simple idea proved extremely difficult and discussions about the Mises' *frequentist theory* never ended. Here is his idea. Toss a coin many times and from time to time calculate the frequency (1.7) of heads. After a sufficiently large $\nu$ it will only change within narrow bounds and at $\nu \to \infty$ it will reach some limiting value. It was this value that Mises called statistical probability (of heads).

Infinitely long trials are impossible, but Mises cited a similar approach in physics and mechanics (for example, velocity at a given moment). He also stated that the sequence of the trials (the *collective*) should be irregular (so that its infinite subsequences should lead to the same probability $\hat{p}$ ).

This condition is too indefinite. How many subsequences ought to be tested before irregularity is confirmed? And is it impossible to select *randomly* an excessively peculiar subsequence? Even these superficial remarks show the great difficulties encountered by the frequentist theory; nevertheless, naturalists have to issue from statistical probability.

Yes, it is theoretically imperfect, although mathematicians came to regard it somewhat milder (Kolmogorov 1963, p. 369). I ought to add that (Uspensky et al 1990, § 1.3.4)

*Until now, it proved impossible to embody Mises' intention in a definition of randomness satisfactory from any point of view*.

**1.1.4.** *Independence of Events and Observations.* Events *A* and *B* are independent if (1.1.1-3)

$P(AB) = P(A)P(B)$,

otherwise

$P(AB) = P(A)P(B/A)$.

Switch *A* and *B*, then

$P(AB) = P(B)P(A/B)$.

A remarkable corollary follows: if *A* does not depend on *B*, then *B* does not depend on *A*; independence is mutual (De Moivre (1718/1756, p. 6):

*Two events are independent, when they have no connection one with the other, and that the happening of one neither forwards nor obstructs the happening of the other.*

*Two events are dependent, when they are so connected together as that the probability of either's happening is altered by the happening of the other.*

The proof of mutuality of independence (already evident in that definition) is simple. According to the condition, $P(A/B) = P(A)$, then by formulas (1.3) and (1.2)

$P(AB) = P(B)P(A)$, $P(B/A) = P(B)$, QED.

Here, however, is a seemingly contradicting example. Suppose that the weather during a summer week in some town is random. Then the random sales of soft drinks there depend on it although there simply cannot be any inverse dependence. But weather and sales cannot be here considered on the same footing.

De Moivre (1711, Introduction) was the first to mention independence, see also just above. Later classics of probability theory mentioned independence of events as well (see below), but some authors forgot about it. The situation had abruptly changed since Markov investigated his *chains* (§ 5.2) and thus added an additional direction to the theory, the study of dependent random events and variables.

Gauss (1823, § 18) stated that if some observation was common to two functions of the results of observations, the errors of these latter will not be independent from each other. He added (for some reason, only in § 19) that those functions were linear. Without this restriction his statement would have contradicted the Student – Fisher theorem about the independence of the sample mean and variance in case of the normal distribution.

Also dependent, as Gauss (1828, § 3) thought, were the results of adjustments. Thus, after the observed angles of a triangle were corrected, and their sum became strictly equal to its theoretical value, these adjusted angles were not anymore independent; they are now somehow connected by their unavoidable residual errors. Note that Gauss had thus considered independence of functions of random variables (§ 1.2.3).

Geodesists invariably (and without citing Gauss) kept to the same definition. Thus, in the Soviet Union separate chains of triangulation had bases and astronomically determined directions on both ends. Therefore, after their preliminary adjustment they were included in a general adjustment as independent entities. True, the bases and those directions were common to at least two chains, but they were measured more precisely than the angles.

Bernstein (1946, p. 47) offered an instructive example showing that pairwise independence of, say, three events, is not sufficient for their mutual independence.

## 1.2. Randomness and Random Variables

**1.2.1.** *Randomness*. In antiquity, randomness was a philosophical notion, then became a mathematical concept as well. Aristotle included it in his doctrine of causes; here are two of his celebrated examples.

1) Digging a hole for a tree, someone finds a buried treasure [not a rusty nail!] (*Metaphysics* 1025a).

2) Two men known to each other meet suddenly (*Physics* 196b30); two independent chains of events *suddenly* intersected.

These examples have a common feature: a small change in the action(s) of those involved led to an essential change: the treasure would have remained buried, there would have been no meeting. Many ancient authors imagined chance just as Aristotle did whereas Cournot (1843, § 40) mentioned the second example anew.

The pattern small change – essential consequences became Poincaré's (1896/1987, pp. 4 – 6) main explanation of randomness, although he specified: when equilibrium is unstable. Here is his or, rather, Cournot's (1843, § 43) example: a right circular cone standing vertically on its vertex falls in a random direction. A similar example is due to Galen (1951, p. 202), a Roman physician and naturalist, $2^{nd}$ century:

*In those who are healthy […] the body does not alter even from extreme causes; but in old men even the smallest causes produce the greatest change.*

*Corruption of nature's aims* was another cause of randomness. Kepler (1618 – 1621/1952, p. 932) established that planets move along elliptical orbits whereas nature, as he thought, aimed at circular orbits. Complete perfection was not attained. Only Newton proved that the ellipticity followed from his law of universal gravitation and that the eccentricity of an orbit was determined by the planet's initial velocity.

Following Kepler and Kant, Laplace (1796/1884, p. 504) somehow concluded that these eccentricities had been caused by variations of temperatures and densities of the diverse parts of the planets.

A mathematical theory cannot however be based on encounters or nature's aims. I leave aside very interesting but occurring much ahead of their time and therefore unsuccessful attempts mathematically to determine randomness (Lambert 1771, §§ 323 – 324; 1772 – 1775), see also Sheynin (1971b, pp. 245 – 246). Modern attempts deal with infinite (and even finite) sequences of zeros and unities such as

0, 0, 1, 1, 1, 1, 1, 0, 1, 1, 0, …

Is it random or not? Such questions proved fundamental. They were approached in various ways, but are far from being solved. For a finite sequence that question is still more complicated. In any case, the beginning of an infinite sequence ought to be irregular so that irregularity (as Mises also thought) is an essential property of randomness.

In philosophy, randomness is opposed to necessity; in natural sciences Poincaré (1896/1912, p. 1) described their dialectic:

*Dans chaque domaine, les lois précises ne décidaient pas de tout, elles traçaient seulement les limites entre lesquelles il était permis au hasard de se mouvoir.*

He did not regrettably mention regularities of mass random events. It is also appropriate to recall the celebrated Laplace's (1814/1995, p. 2) statement allegedly proving that he rejected randomness:

*An intelligence that, at a given instant, could comprehend all the forces by which nature is animated [...], if, moreover, it were vast enough to submit these data to analysis, would encompass [...] the movements of the greatest bodies and those of the slightest atoms. [...] Nothing would be uncertain, and the future, like the past, would be open to its eyes.*

Such intelligence is impossible. Then, there exist unstable motions, responding to small errors of the initial conditions (see above) and perhaps half a century ago a mighty generalization of the former phenomenon, the chaotic motion, was discovered and acknowledged. Finally, Maupertuis (1756, p. 300) and Boscovich (1758, § 385) kept to the same *Laplacean determinism*.

*Allegedly proving ...* Perhaps Laplace's entire astronomical investigations and certainly all his statistical work refute his statements (which really took place) denying randomness.

**1.2.2.** *Cause or Chance*? What should we think if a coin falls on the same side 10 or 20 times in succession? Common sense will tell us: the coin was imperfect. Nevertheless, we will discuss this example. Indeed, after the appearance, in mid-19$^{th}$ century, of the non-Euclidean geometry we may only trust common sense in the first approximation.

Denote heads and tails by + and –. After two tosses the outcomes can be + +, + –, – + and – –, all of them equally probable. After the third toss the outcome + + becomes either + + +, or + + –. In other words, the outcome + + + is not less probable than any of the other 7, and it is easy to see that a similar conclusion remains valid at any number of tosses. Of course 10 heads in succession are unlikely, but all the other possible outcomes will be just as unlikely.

So let us refer to Laplace (1776, p. 152; 1814/1995, p. 9), who discussed the so-called D'Alembert – Laplace problem:

*Suppose we laid out [...] the printer's letters <u>Constantinople</u> in this order. We believe that this arrangement is not due to chance, not because it is less possible than other arrangements. [...] [S]ince we use this word it is incomparably more probable that someone has arranged the preceding letters in this order than that this arrangement happened by chance.*

No formulas can help us and Laplace had to turn to common sense. In our example, we may conclude that someone had done something so that the coin always falls on the same side. Common sense did not let us down. In 1776, Laplace selected the word *Infinitesimal*; it was D'Alembert (1767, pp. 245 – 255) who wrote *Constantinople*. His considerations were not as reasonable.

In general, the *cause or chance* problem compels us to separate somehow equally possible cases (if they exist) into ordinary and remarkable; *Constantinople* was indeed a remarkable arrangement. Kepler was an astrologer as well (and called himself the founder of an allegedly scientific astrology which only admitted a correlative influence of the *stars* on human beings). He (1601, § 40/1979, p. 97) added three aspects (remarkable mutual positions of the heavenly bodies) to the five recognized by the ancients and he (1604/1977, p. 337) also was *not willing to ascribe* the appearance of a New star *to blind chance* [...] and considered it *a great wonder*.

Another related subject is the superstition and self-delusion peculiar to gamblers (and not only to them). A ball is rolled along a roulette wheel and stops with equal probability at any of the 37 positions 0, 1, ..., 35, 36. Gamblers attempt to guess where exactly will the ball stop and the winner gets all the stakes; however, if the ball stops at 0, the stakes go the banker. This is the simplest version of the game.

Now suppose that the ball stopped at 18 three times in succession; should a gambler take this fact into account (and how exactly)?

Petty (1662/1899, vol. 1, p. 64) resolutely opposed games in chance (considered that playing as such was a superstition): *A lottery [...] is properly a tax upon unfortunate self-conceited fools*. Montmort (1708/1980, p. 6) and other authors noted the gamblers' superstitions; and here is Laplace (1814/1995, p. 92) commenting on a similar event:

*When one number has not been drawn for a long time* [...], *the mob is eager to bet on it*.

But it was Bertrand (1888, p. XXII) who dealt the final blow (although did not convince the gamblers): *Elle* [the roulette] *n'a ni conscience ni mémoire. Play, but do not retrieve your losses* (a Russian saying)! It means: play if you cannot abstain from gambling, but never risk much. Arnauld & Nicole (1662/1992, p. 332) warned against expecting large gains (and risking much!).

Laplace (Ibidem, p. 93) also mentioned the general public' superstitions:

*I have seen men, ardently longing for a son* [...]. *They fancied that the boys already born* [during a certain month] *made it more probable that girls would be born next*.

Finally, I note that Laplace (p. 93) saw no advantage in repeatedly staking on the same number. This brings us to martingales, but I will not go thus far.

**1.2.3.** *Random Variable.* This is the central notion of the theory of probability. Here is the simplest definition of a discrete random variable: *A variable taking various discrete values, each with some probability*. Denote these values by $x_1, x_2, ..., x_n$. The sum of their probabilities $p_1, p_2, ..., p_n$ should be unity. Considered together, those values and probabilities are the random variable's *law of distribution*. A random event can be understood as a random variable having $n = 2$.

The case of $n \to \infty$ is also possible; it can be realized in the discrete way, for example, if $x_1 = 1, x_2 = 2, x_3 = 3, ...$, with a *countable* number of the values, or, if a continuous random variable is considered, that number is uncountable. Example: all the uncountable values in interval [0; 1]. A new circumstance appears when there are infinitely many values: an event having a zero probability is possible. Indeed, select any point, say, in the appropriate interval. The probability of choosing any given point is certainly zero, but we did select some point! The geometric probability (§ 1.1.2) can be recalled here.

*A random variable* (or its generalization, which we will not discuss) *or a random event ought to be present in each problem of the theory of probability*. Thus, the outcome of a dice-fall is a random variable; it takes 6 values, each with its own probability (here, they are identical).

Many interesting examples of random variables can be provided. Thus, in the beginning of the 17th century the participants in the celebrated Genoese lottery could guess 1, 2, …, 5 numbers out of 90. The gains increased with those numbers, but the more did the gambler hope for, the heavier was he punished (his expected gain rapidly decreased). This did not follow from any mathematical theorem, but was caused by the organizers' greed.

The random variable involved (the random gain) had 5 values with definite probabilities although only a handful of people had been able to calculate them. Then, from 1662 onward (Graunt), human lifespan began to be studied. In 1756 and 1757 Simpson effectively introduced random variables into the future theory of errors and until about the 1930s this new direction of research had remained the main subject of probability theory. Simpson assumed that the chances of the (random) errors corrupting each measurement (of a given series) are represented by some numbers; the result of measurement thus became a possible value of some random variable and a similar statement held for all of them taken together.

A formal introduction of the random variable was due to Poisson (1837, pp. 140 – 141 and 254) who still called it by a certainly provisional term *chose A*. The proper term, random variable, did not come into general use all at once. Perhaps its last opponent was Markov (letter to Chuprov of 1912; Ondar 1977/1981, p. 65):

*Everywhere possible, I exclude the completely undefined expression <u>random</u> and <u>at random</u>. Where it is necessary to use them, I introduce an explanation corresponding to the pertinent case.*

He had not however devised anything better and often wrote *indefinite magnitude*, which was hardly better. Markov had not applied the terms *normal distribution* or *correlation coefficient* either!

In a certain sense, the entire development of the theory of probability consisted in an ever more general understanding of *random variable*. At first, randomness in general had been studied (actually, a random variable with a uniform law of distribution, see § 2.2.1) as contrary to necessity, then random variables having various distributions, dependent variables and random functions, cf. § 5.1. The level of abstraction in the theory gradually heightened (the same is true about the development of mathematics in general). It is well known that, the higher became that level (i. e., the further mathematics moved away from nature), the more useful it was. Complex numbers and functions of complex variables are absolutely alien to nature, but how useful they are in mathematics and its applications!

## Chapter 2. Laws of Distribution of Random Variables, Their Characteristics and Parameters
### 2.1. Distribution Function and Density

For describing a continuous random variable (call it ξ) we need to determine its law of distribution as it was done in § 1.2.3 for discrete variables. Denote by $F(x)$ the probability of its being less than some $x$:

$P(ξ < x) = F(x).$

This $F(x)$ is called the distribution (integral) function of ξ. If ξ takes any value from $-∞$ to $∞$, then

$P(ξ < -∞) = F(-∞) = 0, P(ξ < ∞) = F(∞) = 1.$

Choose now two arbitrary points, $x_1$ and $x_2$, $x_2 > x_1$, then

$P(ξ < x_2) ≥ P(ξ < x_1)$ or $F(x_2) ≥ F(x_1).$

Indeed, $P(-∞ < ξ < x_2)$ cannot be lower than $P(-∞ < ξ < x_1)$. And if a random variable takes no values on interval $[x_1; x_2]$ (but remains continuous beyond it), then

$P(ξ < x_2) = P(ξ < x_1)$ or $F(x_2) = F(x_1).$  (2.1)

And so, in any case, the function $F(x)$ does not decrease and if (2.1) does not take place, it increases. Note also that

$F(x_2) - F(x_1) = P(ξ < x_2) - P(ξ < x_1).$  (2.2)

Integral distribution functions began to be applied in the $20^{th}$ century, although they fleetingly appeared even in 1669. Pursuing a methodical aim, Huygens (1669/1895, between pp. 530 и 531) drew a graph of a function whose equation can be written as

$y = 1 - F(x), 0 ≤ x ≤ 100.$

The curve described the human lifespan (ξ), the probability of $P(ξ ≥ x)$, but it was not based on numerical data. In 1725, De Moivre studied the same probability, and similarly Clausius (1858/1867, p. 268) investigated the probability of the free path of a molecule to be not less than $x$.

Until distribution functions really entered probability, continuous random variables had been described by *densities* φ(x) of their distributions (of their probability). Consider an infinitely short interval $[x_1; x_1 + dx_1]$. A random variable takes there a value depending on $x_1$; we may say, takes one and the same value φ($x_1$). On the adjacent interval of the same length on the right side the value of that variable may be assumed equal to φ($x_2$), $x_2 = x_1 + dx_1$. Thus we get a series of values φ($x_1$), φ($x_2$), … and can describe the relation of this function, φ(x), the density, with $F(x)$:

$$F(x_n) = \int_{-∞}^{x_n} φ(x)dx, \quad F(x_1) = \int_{-∞}^{x_1} φ(x)dx, \quad F(x_n) - F(x_1) = \int_{x_1}^{x_n} φ(x)dx.$$

These formulas additionally explain equality (2.2). Strictly speaking, by definition,

$F'(x) = \varphi(x),$

but the essence of φ(x) as stated above certainly holds. In more simple examples the density is a continuous function existing on a finite or infinite interval; according to its definition, the area *under* the density curve is unity.

Under, above, to the left or to the right are non-mathematical expressions, but we will apply them even without italics.

Instead of random variables themselves the theory of probability studies their distribution functions or densities just as trinomials

$f(x) = ax^2 + bx + c, a \neq 0$

are studied in algebra. Given the parameters *a*, *b* and *c*, we can determine whether the roots of the trinomial are real (coinciding or not) or complex, can draw its graph. The same way we determine the behaviour of random variables. But where are the parameters of densities or distribution functions?

Consider a function *f(x)*. We may write it down as *f(x; a; b; c)* and thus show that its argument is the variable *x*, but that its behaviour is also determined by parameters constant for each given function (for each trinomial). The density and the distribution function also have parameters peculiar to each random variable. As a rule, statisticians estimate those parameters. Suppose that we have a continuous triangular distribution (assumptions of such kind should be justified) with an unknown parameter *a* (see § 2.2.2). It is required to *estimate* it, to establish for it some (sample) value $\hat{a}$, which is only possible when having appropriate observations of the random variable, and to determine the possible error of that estimate. If there are two parameters, certainly both should be estimated.

## 2.2. Some Distributions

**2.2.1.** *The uniform distribution*. A random variable having this distribution takes all its values with the same probability. Thus, all the 6 outcomes of a die-fall are equally probable. A continuous random variable takes identical values on some interval. The area under this interval should be unity; for interval [– *a, a*] the density will therefore be

$\varphi(x) = 1/a = \text{Const}$

and *a* can be considered the parameter of this distribution.

**2.2.2.** *The continuous triangular distribution* is usually even. So let it cut the *x*-axis at points A (– *a*, 0) and C (*a*, 0). The density is the broken line ABC with AB and BC being the equal lateral sides of the isosceles triangle ABC. The area under it is unity, so we have B(0, 1/*a*).

The only parameter of this distribution is obviously *a* since only it determines the coordinates of all the points A, B and C. I described the triangular distribution mostly since it was easy to establish the meaning of its parameter. It was introduced by Simpson (§ 1.2.3).

**2.2.3.** *The binomial distribution.* We all remember the formula of the square of the sum of two numbers and some of us even managed to remember the formula for the cube of the same sum. However, there exists a general formula for natural exponents $n = 1, 2, 3, \ldots$:

$$(p + q)^n = p^n + C_n^1 p^{n-1} q + C_n^2 p^{n-2} q^2 + \ldots + C_n^{n-1} pq^{n-1} + q^n. \quad (2.3)$$

We are only interested in the particular case of $p + q = 1$, that is, in those magnitudes which describe the probabilities of contrary events. Here, $C_n^k$ is the number of combinations of $n$ taken $k$ at a time:

$$C_n^k = \frac{n(n-1) \ldots (n-k+1)}{k!}, \quad C_n^k = C_n^{n-k}.$$

The numerator has the same number of multipliers as the denominator. Thus,

$$C_5^3 = \frac{5 \cdot 4 \cdot 3}{3!}, \quad 3! = 1 \cdot 2 \cdot 3.$$

Required now is the probability of casting a unity twice when rolling four dice (or rolling one die four times). Cast a die once, and the probability of a unity is $p = 1/6$, whereas the probability of all the other outcomes is $q = 5/6$. And now consider a binomial $[(1/6) + (5/6)]$ raised to the fourth power:

$$[(1/6) + (5/6)]^4 = [1/6^4](1 + 4 \cdot 1^3 \cdot 5 + 6 \cdot 1^2 \cdot 5^2 + 4 \cdot 1 \cdot 5^3 + 5^4).$$

The term $6 \cdot 1^2 \cdot 5^2$ will correspond to the probability sought since it, and only it, includes the multiplier $1^2$, denoting the studied outcome (and another outcome). That probability is $6[1/6^4] \cdot 1^2 \cdot 5^2 = 25/6^3 = 25/216$. We have thus taken into account the number (6) of the possible successive casts (the number of combinations of 4 elements taken 2 at a time). Neglecting this coefficient 6, we would have obtained the probability sought when the successful casts were fixed; for example, if the unity should have occurred in the first and the third roll.

The number of trials $n$ and the ratio $p/q$ can be chosen as the parameters of the binomial distribution (2.3). It is not necessary to choose both $p$ and $q$ since only one of these magnitudes is independent ($p + q = 1$). The example above shows that each term of the binomial expansion (2.3) is the probability

$$p(x) = C_n^k p^{n-k} q^k, \quad x = 0, 1, 2, \ldots, n$$

that the studied random event will occur $k$ times in whichever $n$ trials. The frequency is also essential, see § 2.4.1.

Interesting examples of the binomial distribution include the studies of the sex ratio at births, cf. § 4.2. Its generalization is the multinomial distribution with each trial concerning a random variable taking several values rather than a random event. It is therefore described by a multinomial

$$(a + b + c + \ldots)^n.$$

Pertinent qualitative reasoning without mentioning probabilities were due to Maimonides (Rabinovitch 1973, c. 74):

*Among contingent things some are very likely, other possibilities are very remote, and yet others are intermediate.*

**2.2.4.** *The normal distribution*. The function

$$\varphi(x) = \frac{1}{\sigma\sqrt{2\pi}} \exp[-\frac{(x-a)^2}{2\sigma^2}], \; -\infty < x < \infty, \quad (2.4)$$

is the density of the normal distribution. The stochastic meaning of the two of its parameters, $a$ and $\sigma > 0$, is described in § 2.4.2. The corresponding distribution function is

$$F(z) = \frac{1}{\sigma\sqrt{2\pi}} \int_{-\infty}^{z} \exp[-\frac{(x-a)^2}{2\sigma^2}]dx.$$

Let $a = 0$ and $\sigma = 1$, then, in the *standard* case,

$$F(z) = \frac{1}{\sqrt{2\pi}} \int_{-\infty}^{z} \exp[-\frac{x^2}{2}]dx. \quad (2.5)$$

It is however more convenient to tabulate the function

$$F(z) = \frac{1}{\sqrt{2\pi}} \int_{0}^{z} \exp[-\frac{x^2}{2}]dx. \quad (2.6)$$

Indeed, the integrand in formula (2.5) is an even function so that the integrals (2.5) within $(-\infty; 0]$ and $[0; +\infty)$, are equal to each other and equal to 1/2; within, say, $(-\infty; -1]$ the integral (2.5) is equal to the difference between 1/2 and integral (2.6) at $z = 1$. The value of the function (2.6) at $z \approx 3$ is already 0.499; if $z \to +\infty$ its value is 1/2, or, which is the same, within infinite limits its value is unity, as it should be.

The utmost importance of the normal distribution follows from the so-called *central limit theorem* (CLT), a term due to Polya (1920):

*The sum of a large number of independent random variables, each of them only to a small degree influencing that sum, is distributed normally.*

It was Pearson, who, in 1893, definitively introduced the term *normal distribution* in order to avoid naming it after Gauss (1809) or

Laplace who extensively applied it after non-rigorously proving several versions of the CLT. Galton applied that term before Pearson, but the first to suggest it was Peirce (1873, p. 206).

De Moivre (§ 4.2) considered the appearance of the normal law from a binomial distribution and thus proved a particular case of the CLT. Many authors not to mention Laplace had proved various versions of the CLT, but its rigorous proof was due to Markov and Liapunov, not even to Chebyshev.

Denote the probabilities of a male and female births by $p$ and $q$ and neglect all the other possible births so that $p + q = 1$. Then the probabilities of some number of male births (or of this number remaining within some bounds) can be calculated by means of the normal distribution. This was indeed De Moivre's immediate aim. From 1711 onward the parameter $p/q$ became an object of numerous studies (§ 4.2).

About 1874 Galton (1877) invented the so-called *quincunx*, a device for visually demonstrating the appearance of the normal distribution as the limiting case of the uniform law. Shot was poured through several (say, 20) lines of pins, and each shot 20 times deviated with the same probability to the right or to the left and finally fell on the floor of the device. Thus appeared a normal curve. A special feature of that device was that it showed that the normal law was stable (§ 6.1).

**2.2.5.** *The Poisson distribution.* The law of this discrete distribution (Poisson 1837, p. 205) can be written down as

$$P(x) = \frac{a^x}{x!} e^{-a}, \quad x = 0, 1, 2, \ldots$$

The sum of the probabilities $P(x)$ over all the infinite set of the values of $x$ is 1, as it should be. Indeed, $e^{-a}$ is the common multiplier and

$$\sum_{x=0}^{\infty} \frac{a^x}{x!} e^{-a} = e^{-a}(1 + \frac{a}{1!} + \frac{a^2}{2!} + \frac{a^3}{3!} + \ldots) = e^{-a} e^a = 1.$$

Here is an interesting pattern leading to the Poisson distribution: points are entered on an interval according to a uniform distribution, one by one, independently from each other. It occurs that the number of points situated on some part of that interval obeys the Poisson distribution. Example: the number of calls entering an exchange. Its functioning can therefore be stochastically studied.

Suppose an exchange serves 300 subscribers and the hourly probability of one of them speaking is $p = 0.01$. What will be the probability of four or more independent calls made during an hour? The conditions for the appearance of the Poisson distribution are met, and $a = pn = 3$. Then

$$P(\xi \geq 4) = \sum_{x=0}^{\infty} \frac{a^x}{x!} e^{-a} - P(\xi = 0) - P(\xi = 1) - P(\xi = 2) - P(\xi = 3).$$

The sum is unity (see above) and the other terms are easily calculated.

Another example: the distribution of the stars over the sky (Michell 1767). If they are distributed uniformly (on a sphere rather than interval), some of them will be very close to each other (double, triple, … stars). Even then many such stars had been known, and Michell questioned whether this occurred randomly or not. What is the probability that two stars out of all of them are situated not more than 1° apart?

Newcomb (1860, pp. 427 – 429) applied the Poisson distribution to derive the probability that some small part of the celestial sphere contains *s* stars out of *n* uniformly distributed across the celestial sphere. In a sense, it is this distribution that best describes a *random* arrangement of many points. Its parameter is obviously *a*.

In 1898 Bortkiewicz introduced his *law of small numbers*, and for a few decades it had been considered as the main law of statistics. Actually, it only popularized the then yet little known Poisson distribution which is what Kolmogorov (1954) stated but did not justify his opinion and I (2008) proved that he was correct. Botkiewicz's contribution is deservedly forgotten although mostly owing to previous more particular criticisms.

**2.2.6.** *The hypergeometric distribution.* It is important for acceptance inspection of mass production, see below. Consider the Additional problem No. 4 (Huygens 1657) first formulated by Pascal. Given, 12 counters, 4 of them white (as though defective). Required is the probability that 3 white counters occur among 7 counters drawn without replacement.

Well, actually the entire batch should be rejected, but nevertheless I go ahead following Jakob Bernoulli (1713, part 3, problem 6), although applying the hypergeometric distribution. Huygens, it ought to be added, provided the answer, but not the solution. Denote the conditions of the problem: $N = 12$, $M = 4$, $n = 7$, $m = 3$. Simple combinatorial reasoning lead to a formula which is indeed the formula of that distribution:

$$P(\xi = m) = C_M^m C_{N-M}^{n-m} \div C_N^n.$$

## 2.3. The Main Characteristics of Distributions

**2.3.1.** *Expectation.* For a discrete random variable $\xi$ it is the sum of the products of all its values $x_1, x_2, …, x_n$ by their probabilities $p_1, p_2, …, p_n$:

$$E\xi = \frac{p_1 x_1 + p_2 x_2 + … + p_n x_n}{p_1 + p_2 + … + p_n}. \tag{2.7}$$

The denominator is naturally unity. Laplace (1812/1886, p. 189) added the adjective *mathematical* to expectation so as to distinguish it from the then topical but now forgotten moral expectation (see below). This adjective is regrettably still applied in French and Russian literature.

Expectation can be considered a natural ersatz of a random variable, as though its mean value; in the theory of errors, it corresponds to the generalized arithmetic mean. Denote observations by $x_1, x_2, \ldots, x_n$, and their weights (worth) by $p_1, p_2, \ldots, p_n$. By definition their mean is

$$\bar{x} = \frac{p_1 x_1 + p_2 x_2 + \ldots + p_n x_n}{p_1 + p_2 + \ldots + p_n}, \qquad (2.8a)$$

although the denominator is not 1 anymore. If all the weights are identical

$$\bar{x} = \frac{x_1 + x_2 + \ldots + x_n}{n}. \qquad (2.8b)$$

In § 2.6 I mentioned the selection of bounds covering a measured constant as practised by ancient astronomers. Here, I note that they did not choose any definite estimator, such as the arithmetic mean; they had applied qualitative considerations and thought about convenience of subsequent calculations. For observations corrupted by large errors this tradition makes sense.

So when had that mean become the standard estimator? While selecting a mean of four observations, Kepler (1609/1992, p. 200/63) chose a generalized mean (2.8a) rather than *the letter of the law*, i. e., as I understand him, rather than the mean (2.8b), see Sheynin (1993b, p. 186).

The mean (2.8a) had sometimes been applied with posterior weights $p_i$, equally decreasing on either side of the middle portion of the observations. This choice is hardly useful since, first, these weights are necessarily subjective; and, second, since that estimator only provided a correction of the mean (2.8a) for the unevenness of the sample density of probability of the observational errors.

The expectation (2.7) and the arithmetic mean (2.8) nevertheless essentially differ. The former is a number since it presumably contains all the values of a random variable, whereas the latter is compiled from the results of observations unavoidably corrupted by random errors (as well as by systematic errors, but now we do not need them) and is therefore a random variable as well, as though a sample value of the unknown expectation. Its error ought to be estimated and *a similar remark will also apply to other characteristics of a random variable*. At the same time the arithmetic mean is assumed as the value of the measured constant (§ 6.2). Note that notation $\bar{x}$ for the values of $x_i$ is standard.

For a continuous random variable the expectation is expressed by the integral

$$E\xi = \int_a^b x\varphi(x)dx. \qquad (2.9)$$

Points $a$ and $b$ are the extreme points of the domain of the density $\varphi(x)$ and possibly $a = -\infty$, and $b = \infty$.

Expectation had begun to be applied before probability was. It first appeared, apparently being based on intuitive and subjective chances and in everyday life rather than in science. Maimonides (Rabinovitch 1973, p. 164): *A marriage settlement* [insurance for a woman against divorce or death of husband] *of 1000 zuz can be sold at a present value of 100*, [but] *if the face value were 100 it could not be sold for 10 but rather for less*. Large (though not more likely) gains had been considered preferable, and the same subjective tendency is existing nowadays (and the organizers of lotteries mercilessly take advantage of it). Similar ideas not quite definite either and again connected with insurance appeared in Europe a few centuries later (Sheynin 1977, pp. 206 – 209).

The theory of probability which *officially* originated in 1654, in the correspondence of Pascal and Fermat, effectively applied expectation. Here is one of their main problems which they solved independently from each other. Gamblers A and B agree to play until one of them scores 5 points (not necessarily in succession) and takes both stakes. For some reason the game is interrupted when the score was 4:3 to A. So how should they share the stakes?

Even then that problem was venerable; there are indications that a certain mathematician had solved it at least in a particular case. Note that sharing the stakes proportionally to 4:3 would have been fair when playing chess, say, i. e., when the gamblers' skill is decisive. In games of chance, however, everything depends on chance and the past cannot influence the future (cf. § 1.2.2).

Here is the solution. Gambler A has probability $p_1 = 1/2$ (Pascal and Fermat kept to chances) of winning the next play; he can also lose it with the same probability but then the score will equalize and the stakes should be equally shared. A's share (the expectation of his gain) will therefore be $1/2 + 1/4 = 3/4$ of both stakes. The expectation of the second gambler is therefore 1/4 of both stakes and it could have been calculated independently.

It was a man about town, De Méré, who turned Pascal's attention to games of chance (Pascal 1654/1998, end of Letter dated 29 July 1654). He was unable to understand why the probability of an appearance of a six in 4 casts of a die was not equal to that of the appearance of two sixes in 24 casts of two dice as it followed from an old approximate rule. Here, however, are those probabilities:

$$P_1 = 1 - (5/6)^4 = 0.518, P_2 = 1 - (35/36)^{24} = 0.492.$$

So De Méré knew that gamblers had noted a difference of probabilities equal to 0.026. Cf. a similar remark made by Galileo (§ 1.1.1).

Huygens (1657) published a treatise on calculations in games of chance. He formally introduced the expectation in order to justify the sharing of stakes and the solutions of other problems. He substantiated the expediency of applying it for estimating a random variable (a random winning) by reasonable considerations.

Jakob Bernoulli (1713, part 1) however suggested a much simpler justification. Here is a quotation from Huygens and Bernoulli's reasoning (his part 1 was a reprint of Huygens complete with important comments).

*Huygens, Proposition 3*. Having *p* chances to get *a* and *q* chances to get *b* and supposing that all these chances are the same, I obtain

$$\frac{pa + qb}{p + q}. \qquad (2.10)$$

Since *p* and *q* are chances rather than probabilities, their sum is not unity as it was in formula (2.7). And here is Bernoulli. Suppose there are (*p* + *q*) gamblers, and each of *p* boxes contains sum *a*, and each of *q* boxes contains *b*. Each gambler takes a box and all together get (*pa* + *qb*). However, they are on the same footing, should receive the same sum, i. e., (2.10).

As stated in § 1.2.1, a mathematical theory cannot be based on boxes or gamblers, and even De Moivre introduced expectation axiomatically, without justifying it. And so it is being introduced nowadays, although Laplace (1814/1886, p. XVIII) just stated that it is *la seule equitable*.

Several centuries of applications have confirmed the significance of the expectation although in 1713 Nikolaus Bernoulli, in a letter to Montmort published by the latter (Montmort 1708/1713, p. 402) devised a game of chance in which it did not help at all.

Gambler A casts a die … However, the die was very soon replaced by a coin. And so, if heads appears at once, B pays A 1 écu; if heads only appears at the second toss, he pays 2 écus, 4 écus if only at the third toss etc. Required is the sum which B ought to receive beforehand.

Now, A gets 1 écu with probability 1/2, 2 écus with probability 1/4, 4 écus with probability 1/8 etc and the expectation of his gain is

$$1 \cdot (1/2) + 2 \cdot (1/4) + 4 \cdot (1/8) + \ldots = (1/2) + (1/2) + (1/2) + \ldots = \infty. \quad (2.11)$$

However, no reasonable man will agree to pay B any considerable sum and hope for a large (much less, an infinite) gain. He will rather decide that heads will first occur not later than at the sixth or seventh toss and that he ought to pay beforehand those 1/2 écus not more than six or seven times; all the rest infinite terms of the series (2.11) will therefore disappear.

Buffon (1777, § 18) reported that 2048 such games resulted in A's mean gain of 4.9 écus and that only in 6 cases they consisted of 9 tosses, of the largest number of them. His was the first statistical study of games of chance. On a much greater scale Dutka (1988) conducted a similar investigation by applying a computer.

This paradox continued to interest mathematicians up to our time, but it was Condorcet (1784, p. 714) who left the most interesting remark: one game, even if infinite, is still only one trial; many games are needed for stochastically considering them. Freudenthal (1951) independently repeated this remark and additionally suggested that before each game the gamblers ought to decide by lot who will pay whom beforehand.

A similar statement about neglecting low probabilities holds for any game of chance (and any circumstance in everyday life). If there are very large gains in a lottery available with an extremely low probability (which the organizers will definitely ensure), they ought to be simply forgotten, neglected just like the infinite tail of the series (2.11).

But then, how low should a neglected probability be? Buffon (1777, § 8), issuing from his mortality table, suggested the value 1/10,000, the probability of a healthy man 56 years old dying within the next 24 hours. What does it mean for the Petersburg game? We have

$1/2^n = 1/10,000$, $2^n = 10,000$, $n\lg 2 = 4$ and $n \approx 13.3$.

Even that is too large: recall Buffon's experiment in which the maximal number of tosses only amounted to 9. This result also means that 1/10,000 was too low; we may often neglect much higher probabilities and, anyway, a single value for a neglected probability valid in any circumstances should not be assigned at all. And some events (the Earth's collision with a large asteroid) should be predicted with a probability much higher than $(1 - 1/10,000)$. It is not however clear how to prevent such global catastrophes.

Reader! Do you think about such probabilities when crossing the road?

While attempting to solve the paradox of the invented game, Daniel Bernoulli (1738) introduced *moral expectation* (but not the term itself). He published his memoir in Petersburg, and thus appeared the name *Petersburg game.* In essence, he thought that the real value of a gambler's gain is the less the greater is his fortune. He applied his novelty to other risky operations and for some decades it had been widely appraised (but not implemented in practice). At the end of the 19th century economists had developed the theory of marginal utility by issuing from moral expectation.

**2.3.1**-**1.** *The properties of the expectation.* **1)** Suppose that $\xi = c$ is constant. Then

$$\mathrm{E}c = \int_a^b c\varphi(x)dx = c\int_a^b \varphi(x)dx = c.$$

*The expectation of a constant is that very constant.*

**2)** The expectation of a random variable $a\xi$ is

$$\mathrm{E}a\xi = a\int_a^b x\varphi(x)dx = a\mathrm{E}\xi.$$

*When multiplying a random variable by a constant its expectation is multiplied by that very constant.*

**3)** Two (or more) random variables $\xi$ and $\eta$ are given; their densities are $\varphi(x)$ and $\eta(y)$ and the expectation of their sum is sought. It is equal to the double integral

$$E(\xi + \eta) = \int_a^b \int_c^d (x+y)\varphi(x)\psi(y)dxdy =$$

$$\int_a^b \int_c^d x\varphi(x)\psi(y)dydx + \int_a^b \int_c^d y\varphi(x)\psi(y)dxdy.$$

Here, $c$ and $d$ are the extreme points of the domain of the second function and $a$ and $b$ have a similar meaning (see above). Notation $\eta(y)$ instead of $\eta(x)$ does not in essence change anything but transformations become clearer.

The first integral can be represented as

$$\int_c^d \psi(y)dy \int_a^b x\varphi(x)dx = E\xi,$$

since the integral with respect to $y$ is unity. Just the same, the second integral is $E\eta$ and therefore

$$E(\xi + \eta) = E\xi + E\eta.$$

*The expectation of a sum of random variables is equal to the sum of the expectations of the terms. A similar statement can be proved about the difference of random variables*: *its expectation is equal to the difference of the expectations of the terms*.

Note however that differences in such theorems (not only in the theory of probability) are usually not mentioned since by definition subtraction means addition of contrary magnitudes; thus, $a - c \equiv a + (-c)$.

**4)** Without proof: *the expectation of a product of two <u>independent</u> random variables equals the product of their expectations*:

$$E(\xi\eta) = E\xi \cdot E\eta.$$

This property is immediately generalized on a larger number of random variables.

All the properties mentioned above also take place for expectations of discrete random variables.

**2.3.2.** *Variance* is the second main notion characterizing distributions of random variables, their scattering. An inscription on a Soviet matchbox stated: *approximately 50 matches*. But suppose that actually one such box contains 30 matches, another one, 70. The mean is indeed 50, but is not the scattering too great? And what does *approximately* really mean?

Suppose that only some values $x_1, x_2, \ldots, x_n$ of a random variable $\xi$ (a sample of size $n$) are/is known. Then the sample variance of $\xi$ is

$$s^2 = \frac{\sum_{i=1}^n (x_i - \bar{x})^2}{n-1}. \qquad (2.12)$$

It is also called *empirical* since the values of $x_i$ are the results of some experiment or trial.

Why function (2.12) is chosen as a measure of scattering, and why its denominator is $(n - 1)$ rather than $n$? I attempt at explaining it, but first I add that the variance (not sample variance) of the same $\xi$, var$\xi$, of a discrete or continuous variables is, respectively,

$$\sum_{i=1}^{n} p_i (x_i - E\xi)^2, \quad \sigma_\xi^2 = \int_a^b (x - E\xi)^2 \varphi(x) dx, \qquad (2.13)$$

where $a$, $b$ and $\varphi(x)$ have the same meaning as in formula (2.9).

It was Gauss (1823) who introduced the variance as a measure of the scattering of observations. Its choice, as he indicated, is more or less arbitrary, but such a measure should be especially sensitive to large errors, i. e. should include $(x - E\xi)$ raised to some natural power (2, 3, …), and remain positive which excludes odd powers. Finally, that measure should be as simple as possible which means the choice of the second power of that binomial. Actually, Gauss (1823, §§ 37 – 38) had to determine only the sample variance and to apply the arithmetic mean instead of the expectation. Below, I will say more about the advantages of the variance.

Suppose that $x_i$, $i = 1, 2, …, n$, are the errors of observation, then the sample variance will be

[$xx$]/$n$

where [$xx$] is Gauss' notation denoting the sum of the squares of the $x_i$. These errors are however unknown, and we have to replace them by the deviations of the observations from their arithmetic mean. Accordingly, as Gauss proved in the sections mentioned above, the sample variance ought to be represented by formula (2.12). He (1821 – 1823/1887, p. 199) remarked that that change was also demanded by the *dignity of science*.

But suppose that a series of observations is corrupted by approximately the same systematic error. Then those formulas will not take it into considerations, will therefore greatly corrupt reality: the scattering will not perhaps be large although the observations deviated from the measured constant. Gauss himself had directly participated in geodetic observations, therefore did not trust his own formulas (because of the unavoidable systematic errors) and measured each angle until becoming satisfied that further work was useless. Extracts from his field records are published in vol. 9 of his *Werke*, pp. 278 – 281.

Not only the sample variance, but a square root of it (not only $s^2$, but $s$) is applied as well. That $s$ is called *standard deviation*, or, in the theory of errors, *mean square error*.

And now we can specify statements similar to *approximately 50 matches in a box*. Carry out a thankless task: count the matches $x_1$, $x_2$, …, $x_{10}$ in 10 boxes, calculate their mean $\bar{x}$ (their sample mean, since the number of such boxes is immense), the deviations $(x_1 - \bar{x})$,

$(x_2 - \bar{x})$, ..., $(x_{10} - \bar{x})$, and finally the sample variance (or standard deviation). A deviation of some $x_i$ from the approximately promised value that exceeds two mean square errors is already serious.

The expectation of a random variable can be infinite, as in the case of the Petersburg game, and the same can happen with the variance.
*Example.* A continuous random variable distributed according to the Cauchy law

$$\varphi(x) = \frac{2}{\pi(1+x^2)}, \quad 0 \leq x < \infty. \tag{2.14}$$

Note that equalities such as $x = \infty$ should be avoided since infinity is not a number but a variable. Also bear in mind that the distribution (2.14) first occurred in Poisson (1824, p. 278).

Now, the variance. It is here

$$\text{var}\xi = \frac{2}{\pi} \int_0^\infty x^2 \cdot \frac{1}{1+x^2} dx = \frac{2}{\pi} \left[ \int_0^\infty 1 \cdot dx - \int_0^\infty \frac{1}{1+x^2} dx \right]$$

The second integral is

$$\text{arctg} x \Big]_0^\infty = \pi/2,$$

but the first does not exist (and the variance is infinite):

$$\int_0^\infty dx = x \Big]_0^\infty \to \infty.$$

The arithmetic mean of observations, if they are so unsatisfactory that their errors obey the Cauchy distribution, is not better than an isolated observation. Indeed, according to formula (2.16) from § 2.3.2-2 the variance of the mean of $n$ observations is $n$ times less than the variance of a single observation, that is, $n$ times less than infinity and is therefore also infinite.

**2.3.2-1.** *A second definition of variance.* Definition (2.13b) can be written as

$$\text{var}\xi = \int_a^b x^2 \varphi(x) dx - 2 \int_a^b x \text{E}\xi \varphi(x) dx + \int_a^b (\text{E}\xi)^2 \varphi(x) dx.$$

Now, $\text{E}\xi$ is constant and can be separated:

$$\text{var}\xi = \int_a^b x^2 \varphi(x) dx - 2\text{E}\xi \int_a^b x\varphi(x) dx + (\text{E}\xi)^2 \int_a^b \varphi(x) dx.$$

By definition, the first integral is $\text{E}\xi^2$, and the second, $\text{E}\xi$. The third integral is unity according to the property of the density. Therefore,

$$\mathrm{var}\xi = E\xi^2 - 2(E\xi)^2 + (E\xi)^2 = E\xi^2 - (E\xi)^2. \qquad (2.15)$$

This formula is usually assumed as the main definition of variance.
**2.3.2-2.** *The properties of density.*
**1)** *The density of a sum of independent random variables.* By the second definition of variance we have

$$\mathrm{var}(\xi + \eta) = E(\xi + \eta)^2 - [E(\xi + \eta)]^2 =$$
$$E\xi^2 + 2E(\xi\eta) + E\eta^2 - [(E\xi)^2 + 2E\xi E\eta + E\eta^2].$$

Then, according to the fourth property of expectation of independent random variables (§ 2.3.1-1),

$$E(\xi\eta) = E\xi \cdot E\eta$$

so that

$$\mathrm{var}(\xi + \eta) = [E\xi^2 - (E\xi)^2] + [E\eta^2 - (E\eta)^2] = \mathrm{var}\xi + \mathrm{var}\eta.$$

*The variance of a sum of independent random variables is equal to the sum of their variances.*
**2)** *Corollary*: Variance of the arithmetic mean. Given observations $x_1, x_2, \ldots, x_n$ and their arithmetic mean (2.8b)

$$\bar{x} = \frac{x_1 + x_2 + \ldots + x_n}{n}$$

is calculated. Formula (2.12) provides the sample variance of observation $x_i$, but now we need the variance of the mean. By the theorems on the variance of the sum of random variables (the results of observation are random!) and on the product of a random variable by a constant (here, it is $1/n$), we obtain at once a simple but important formula

$$\mathrm{var}\,\bar{x} = \frac{(x_1 - \bar{x})^2 + (x_2 - \bar{x})^2 + \ldots + (x_n - \bar{x})^2}{n(n-1)}. \qquad (2.16)$$

*The variance of the arithmetic mean of n observations is n times less than the variance of each of them.*
Here, like in formula (2.12), we certainly assume that the observations are possible values of one and the same random variable.
**3)** *The variance of a linear function of a random variable.* Suppose that $\eta = a + b\xi$ is a linear function of random variable $\xi$ (and therefore random as well just like any function depending on a random variable). The variance of $\xi$, $\mathrm{var}\xi$, is known and required is $\mathrm{var}\eta$. Such problems occur often enough.
By formula (2.15)

$$\mathrm{var}\eta = E\eta^2 - (E\eta)^2 = E(a + b\xi)^2 - [E(a + b\xi)]^2.$$

The first term is

$$E(a^2 + 2ab\xi + b^2\xi^2) = a^2 + 2abE\xi + b^2E\xi^2.$$

The second term is

$$[Ea + E(b\xi)]^2 = (Ea)^2 + 2EaE(b\xi) + (Eb\xi)^2 = a^2 + 2abE\xi + b^2(E\xi)^2$$

and their difference is $b^2[E\xi^2 - (E\xi)^2]$.

According to formula (2.15) $\mathrm{var}\eta = b^2\mathrm{var}\,\xi$.

And so, an addition of a constant to a random variable does not change the variance, and, when multiplying such a variable by a constant coefficient, its variance is multiplied by the square of that constant:

$$\mathrm{var}(a + \xi) = \mathrm{var}\xi, \quad \mathrm{var}(b\xi) = b^2\mathrm{var}\xi.$$

### 2.4. Parameters of Some Distributions

In § 2.2 we have determined the parameters of a few distributions, but the binomial and the normal laws are still left.

**2.4.1.** *The binomial distribution.* Suppose that $\mu_k$ is a random number of the occurrences of an event in the $k$-th trial, 0 or 1. If the probability of its happening is $p$, then

$$E\mu_k = 1\cdot p + 0\cdot q = p.$$

In a series of $n$ trials that event occurs

$$(\mu_1 + \mu_2 + \ldots + \mu_n) = \mu \text{ times}, \quad E\mu = E\mu_1 + E\mu_2 + \ldots + E\mu_n = pn.$$

Then, see formula (2.15),

$$\mathrm{var}\mu_k = E\mu_k^2 - (E\mu_k)^2.$$

However, $\mu_k^2$ takes the same values, 0 and 1, as $\mu_k$, and with the same probabilities, $p$ and $q$, so that

$$\mathrm{var}\mu_k = p - p^2 = p(1 - p) = pq, \quad \mathrm{var}\mu = \mathrm{var}\mu_1 + \mathrm{var}\mu_2 + \ldots + \mathrm{var}\mu_n = pqn.$$

The magnitudes $E\mu$ and $\mathrm{var}\mu$ characterize the frequency $\mu$. Recall that in § 2.2.3 we discussed the parameters of the binomial distribution proper.

**2.4.2.** *The normal distribution.* It follows from formula (2.4) that the form of the normal curve depends on the value of $\sigma$; the less it is, the more is the area under that curve concentrated in its central part. The values of the random variable $\xi$ close to the abscissa of the curve's maximum become more probable, the random variable as though shrinks.

At $a = 0$ the graph of the density of the normal distribution becomes symmetrical with respect to the $y$-axis so that $a$ is the *location*

*parameter*. Note that this term is applied to any densities whose formula contains the difference $x - a$.

The analytical meaning of both parameters is very simple:

$$a = E\xi, \quad \sigma^2 = \text{var}\xi. \qquad (2.17a, 2.17b)$$

We will prove (2.17a) and outline the proof of (2.17b). We have

$$E\xi = \frac{1}{\sigma\sqrt{2\pi}} \int_{-\infty}^{\infty} x \exp[-\frac{(x-a)^2}{2\sigma^2}]dx.$$

Now, $x = [(x - a) + a]$ and the integral can be written as

$$\frac{1}{\sigma\sqrt{2\pi}} \{\int_{-\infty}^{\infty} (x-a)\exp[-\frac{(x-a)^2}{2\sigma^2}]dx + a\int_{-\infty}^{\infty} \exp[-\frac{(x-a)^2}{2\sigma^2}]dx\}.$$

In the first integral, the integrand is an odd function of $(x - a)$, which, just as $x$, changes unboundedly from $-\infty$ to $\infty$. This integral therefore disappears (the *negative* area under the $x$-axis located to the left of the $y$-axis is equal to the positive area above the $x$-axis located to the right of the $y$-axis).

Then, in the second integral, let

$$\frac{x-a}{\sigma\sqrt{2}} = z, \quad dx = \sigma\sqrt{2}dz, \qquad (2.18)$$

so that it is equal to

$$\int_{-\infty}^{\infty} \exp(-z^2)dz \cdot \sigma\sqrt{2}.$$

Euler was the first to calculate it; without the multiplier $\sigma\sqrt{2}$ it is equal to $\sqrt{\pi}$. Finally, taking into account all three multipliers, $a$, $\sigma\sqrt{2}$ и $1/\sigma\sqrt{2\pi}$, we arrive at $a$, QED.

Now the formula (2.17b):

$$\text{var}\xi = E(\xi - E\xi)^2 = \frac{1}{\sigma\sqrt{2\pi}} \int_{-\infty}^{\infty} (x-a)^2 \exp[-\frac{(x-a)^2}{2\sigma^2}]dx.$$

We have applied here the just derived formula (2.17a). Now we ought to introduce a new variable, see (2.18), and integrate by parts.

### 2.5. Other Characteristics of Distributions

**2.5.1.** *Those replacing expectation.* For a sample (sometimes the only possibility) those characteristics replace the arithmetic mean or estimate the location of the measured constant in some other way.

**2.5.1-1.** *The median.* Arrange the observations $x_1, x_2, \ldots, x_n$ of a random variable in an ascending order and suppose that the thus ordered sequence is $x_1 \leq x_2 \leq \ldots \leq x_n$. Its median is the middlemost observation, quite definite for odd values of *n*. Suppose that $n = 7$, the median will then be $x_4$. For even values of *n* the median will be the halfsum of the two middle terms; thus, for $n = 12$, the halfsum of $x_6$ and $x_7$.

For continuous random variables with density $\varphi(x)$ the median is point $x_0$ which divides the area under the density curve into equal parts:

$$\int_a^{x_0} \varphi(x)dx = \int_{x_0}^b \varphi(x)dx = 1/2.$$

In other words, the median corresponds to equality $F(x) = 1/2$. To recall: the entire area under the density curve is unity; *a* and *b* are the extreme points of the domain of $\varphi(x)$.

For some densities, as also when the density is unknown, the median characterizes a random variable more reliably then the arithmetic mean. The same is true if the extreme observations possibly are essentially erroneous. Indeed, they can considerably displace the mean but the median will be less influenced.

Mendeleev (1877/1949, p. 156), who was not only a chemist, but an outstanding metrologist, mistakenly thought that, when the density remained unknown, the arithmetic mean ought to be chosen.

Continuous distributions are also characterized by quantiles which correspond to some probabilities *p*, that is, points $x = x_p$ for which $F(x) = p$, so that the median is a quantile corresponding to $p = 1/2$. Its exact location can be not quite certain, cf. the case of the median.

**2.5.1-2.** *The mode.* This is the point (or these are the points) of maximal density. It (one of them) can coincide with the arithmetic mean. Accordingly, the density is called unimodal, bimodal, … or even antimodal (when a density has a point of minimum). In case of discrete random variables the mode is rarely applied.

**2.5.1-3.** *The semi-range* (*mid-range*). This is a very simple but unreliable measure since the extreme values can be considerably erroneous and no other observations are taken into account. It had been widely applied in the 18[th] century for estimating mean monthly values of meteorological elements (for example, air temperatures). It was certainly easier to calculate the mid-range than the mean of 30 values. Interestingly, Daniel Bernoulli (1778, § 10) indicated that he had found it to be *less often wrong* than [he] thought …

**2.5.2.** *Characteristics replacing the variance*

**2.5.2-1.** *The range.* The (sample) range is the difference between the maximal and the minimal measured values of a random variable, cf. § 2.5.1-3. The not necessarily equal differences $(x_n - \bar{x})$ and $(\bar{x} - x_1)$ are also sometimes applied. All these differences are unreliable. In addition to the remarks in that subsection I note that they can well increase with the number of observations; there can appear a value less than $x_1$ or larger than $x_n$.

It is certainly possible to apply instead the fractions $(x_n - x_1)/n$, $(x_n - \bar{x})/n$ and $(\bar{x} - x_1)/n$. The denominator coincides with the possibly increasing number of observations but the numerator changes uncertainly. All the measures mentioned here concern a series of observations rather than a single result.

**2.5.2-2.** *The mean absolute error.* It, just as the probable error (see 2.5.2-3), characterizes a single observation. Denote observations by $x_1$, $x_2$, …, $x_n$, then the mean absolute error will be

$$\sum_{i=1}^{n} |x_i| \div n.$$

It had been applied, although not widely, when treating observations.

**2.5.2-3.** *The probable error.* It was formally introduced by Bessel (1816, pp. 141 – 142) as a measure of precision, but even Huygens (1669/1895), in a letter to his brother dated 28 Nov. 1669, mentioned the idea of a probable value of a random variable. Discussing the random duration of human life, he explained the difference between the expected interval (the mean value derived from data on many people) and the age *to which a person with equal probabilities can live or not*.

Both durations of life should be calculated separately for men and women, which in those times was not recognized. Women *generally* live longer and this possibly compensates them for a life more difficult both in the biological and social sense but they seem to recall this circumstance rather rarely.

Bessel had indeed applied that same idea, repeatedly found in population statistics and, for example, when investigating minor changes in the period of the swings of a pendulum (Daniel Bernoulli 1780). According to Bessel, a probable error of an observation is such that *with equally probability will be either less or larger than the really made error.*

For symmetric distributions the probable error is numerically equal to the distance between the median and the qauntile corresponding to $p = 1/4$ or $3/4$; it is the probability that an observation thus deviates in either side from the median. For the normal distribution that distance is $0.6745\sigma$, and many authors had understood (still understand?) that relation as a universal formula or had tacitly thought that they have dealt with the normal distribution.

Moreover, I am not sure that there exists a generally accepted definition of the probable error suitable for asymmetric distributions, i. e., when the distances from the median to the quantiles $p = 1/4$ and $3/4$ do not coincide. If in such cases the probable error is still meaningful, it is perhaps permissible to say that it is equal to half the distance between those quantiles.

The idea of the probable error is so natural that that measure became universally adopted whereas, perhaps until the second half of the 20[th] century, the mean square error had been all but forgotten. In the third (!) edition of his serious geodetic treatise Bomford (1971, pp. 610 – 611) *reluctantly* abandoned it and went over to the mean square error.

So why is the latter better? We may bear in mind that the probable error is connected with the median which is not always preferable to the arithmetic mean. Then, it, the mean square error (or, rather, the variance), is the most reliable measure. The variance (we may only discuss the sample variance) is a random variable, it therefore has its own sample variance. True, as mentioned above, a similar remark is applicable to any sample measure (in particular, to the arithmetic mean). However, unlike other measures of scattering, the variance of the variance is known, first derived by Gauss (1823, § 40). True, he made an elementary mistake corrected by Helmert (1904), then independently by Kolmogorov et al (1947).

One circumstance ought to be however indicated. Practically applied is not the variance, but its square root, the standard deviation (the mean square error); and if the variance of the variance is $a$, it does not at all mean that the variance of the latter is $\sqrt{a}$; for that matter, it is only known for the normal distribution, see below. Again, the sample variance is an *unbiased* estimate of the general, of the *population* variance which means that its expectation is equal to that variance, whereas the sample standard deviation has no similar property. Recall that Gauss (§ 2.3.2) remarked that the formula for the sample variance had to be changed; now I additionally state that he had thus emphasised the essential role of unbiasedness although currently it is much less positively estimated.

**2.5.2-4.** *An indefinite indication of scattering.* We sometimes meet with indications such as *This magnitude is equal to a ± c*. It can be understood as … *equal to any value between a – c and a + c*, but it is also possible that $c$ is not the maximal but, for example, the probable error. And, how was that $c$ obtained? We have approached here the important subject of interval estimation.

### 2.6. Interval Estimation

Denote some parameter of a function or density by $\lambda$ and suppose that its sample value $\hat{\lambda}$ is obtained. Required is an estimate of the difference $|\hat{\lambda} - \lambda|$. Its *interval* estimation means that, with α and δ being indicated,

$$P(|\hat{\lambda} - \lambda| < \delta) > 1 - \alpha.$$

Now, we may state that the *confidence interval* $[\hat{\lambda} - \delta; \hat{\lambda} + \delta]$ covers the unknown $\lambda$ with *confidence probability* (*confidence coefficient*) $(1 - \alpha)$. This method of estimation is reasonable if α is small (for example, 0.01 or 0.05, but certainly much larger than the Buffon value 1/10,000), and such that δ is also sufficiently small. Otherwise the interval estimation will show that either the number of observations was too small or that they were not sufficiently precise. Note also that in any case other observations can lead to other values of $\hat{\lambda}$ and δ.

Suppose that a constant $A$ is determined by observations. Then, adopting simplest assumptions (Bervi 1899), we may assume that the obtained range $[x_1; x_n]$ covers it with probability

$$P(x_1 \leq A \leq x_n) = 1 - 1/2^{n-1}.$$

I indicated the deficiency of this trick in § 2.5.1-3. Similar conclusions were made by astronomers in the antiquity (Sheynin 1993b, § 2.1).

Issuing from all the existing observations (not only his own) the astronomer selected some bounds (*a* and *b*) and stated that $a \leq A \leq b$. Probabilities had not been mentioned but the conclusion made was considered almost certain.

When determining a constant, any measure of scatter may be interpreted as tantamount to a confidence characteristic. Indeed, suppose that the arithmetic mean $\bar{x}$ of observations is calculated and its mean square error *m* determined. Then the probability $P(\bar{x} - m \leq \bar{x} \leq \bar{x} + m)$ can be established by statistical tables of the pertinent law of distribution as $P(0 \leq \bar{x} \leq \bar{x} + m) - P(0 \leq \bar{x} \leq \bar{x} - m)$; the difference between strict and non-strict inequalities can be neglected. So exactly that *P* is indeed the confidence probability and [$\bar{x} - m$; $\bar{x} + m$], the confidence interval.

## 2.7. The Moments of a Random Variable

This subject can be quite properly included in § 2.6, but it deserves a separate discussion. Moments characterise the density and can sometimes establish it.

The initial moment of order *s* of a discrete or continuous random variable $\xi$ is, respectively,

$$\alpha_s(\xi) = \sum_x x^s p(x), \quad v_s = \int x^s \varphi(x) dx. \qquad (2.19)$$

In the first case, the summing is extended over all the values of *x* having probabilities $p(x)$ whereas the integral is taken within the extreme points of the domain of the known or unknown density $\varphi(x)$ of the continuous random variable.

Also applied are the central moments

$$\mu_s(\xi) = \sum_i (x_i - E\xi)^s p(x_i), \quad \mu_s(\xi) = \int (x - E\xi)^s \varphi(x) dx. \quad (2.20)$$

Both these formulas (the integral is taken between appropriate bounds) can be represented as

$$\mu_s(\xi) = E(\xi - E\xi)^s.$$

Sample (empirical) initial moments for both discrete and continuous random variables certainly coincide:

$$m_s(\xi) = \sum_i x_i^s \div n, \qquad (2.21)$$

where *n* is the number of measured (observed) values of $\xi$.

The central sample moments are

$$m_s(\xi - \bar{x}) = \frac{1}{n} \sum (x_i - \bar{x})^s.$$

The measured (observed) values are often combined within certain intervals or categories. Thus (Smirov & Dunin-Barkovski 1959/1969,

§ 1 in Chapter 3), 70 samples containing 5 manufactured articles each were selected for checking the size of such articles. In 55 samples the size of each of the 5 articles was standard, in 12 of them 2 were non-standard, and in 3, 1 was non-standard:

| Number of samples | 55 | 12 | 3 |
| Number of defective articles | 0 | 1 | 2 |
| Frequencies of the various outcomes | 0.786 | 0.171 | 0.043 |

Here the frequency, for example in the fist column is 55/70.

Many definitions of mathematical statistics have been offered, but only once were statistical data mentioned (Kolmogorov & Prokhorov 1974/1977, p. 721): they denote *information about the number of objects which possess certain attributes in some more or less general set*.

Those numbers above are indeed statistical data; they were separated into sets with differing numbers of defective articles in the samples. Such separation can often be made in several ways; however, if the range of the values of the random variable (the number of defective articles) is sufficiently wide (here, we have a very small range from 0 to 5, but actually even from 0 to 2), there should not be too few sets or intervals. On the other hand, there should not be too many of them either: too many subdivisions of the data is a *charlatanisme scientifique* (Quetelet 1846, p. 278)

And so, when combining the data, formula (2.21) becomes

$$m_s(\xi) = \frac{1}{n} \sum_i n_{x_i} x_i^s, \qquad (2.22)$$

where $n_x$ is the number of the values of the random variable in interval $x$. In our example

$$m_s(\xi) = \frac{1}{70}(55 \cdot 0^s + 12 \cdot 1^s + 3 \cdot 2^s) = 0.786 \cdot 0^s + 0.171 \cdot 1^s + 0.043 \cdot 2^s.$$

The cases of $s = 1$ and 2 (see the very beginning of this section) coincide with expectation and variance respectively and formulas (2.20) correspond to formulas (2.13). *The first moment is the expectation, the second moment is the variance*. But is it possible and necessary to establish something about the other almost infinitely many moments? Suffice it to consider the next two of them.

Suppose that the density $\varphi(x)$ is symmetrical with respect to the $y$-axis. Then for odd values of $s$ the moments

$$v_s = \int_{-\infty}^{\infty} x^s \varphi(x) dx = 0.$$

Indeed, in this case the integrand is the product of an odd and an even function and is therefore odd and the integral is taken between symmetrical bounds.

If some odd moment differs from zero, the density cannot be symmetric (i. e., even) and this moment will therefore characterize the deviation of φ(x) from symmetry. But which moment should we choose as the measure of asymmetry?

All sample moments depend on the observed values of the appropriate random variable, are therefore random variables as well and possess a variance. It is also known that *the variances of the moments of higher orders are larger than those of the first few*. The moments of higher orders are therefore *unreliable*.

It is thus natural to choose the third moment as the measure of asymmetry of the density φ(x) of a random variable; more precisely, the third sample moment since apart from observations we have nothing to go on:

$$m_3(\xi) = \frac{1}{n-1} \sum (x_i - \bar{x})^3. \qquad (2.23)$$

One more circumstance. The dimensionality of the third moment is equal to the cube of the dimensionality of $(x_i - \bar{x})$. For obtaining a dimensionless measure, (2.23) should be divided by $s^3$, see formula (2.12). The final measure of asymmetry of φ(x) is thus

$$s_k = m_3 \div s^3.$$

When discussing that formula (2.12), we indicated why its denominator should be $(n-1)$ rather than $n$. The same cause compelled us to change the denominator in formula (2.23).

Now the fourth moment. For a normal random variable it is $3\sigma^4$, whereas the second moment is $\sigma^2$, see § 2.4.2. For that distribution we therefore have

$$\nu_4/\sigma^4 = 3.$$

If we now calculate the so-called excess (more precisely, the sample excess)

$$\varepsilon_k = m^4/s^4 - 3,$$

its deviation from 0 can be chosen as a measure of the deviation of un unknown density of distribution from the normal law (for which the excess disappears). The excess is here useful since in one or another sense the normal distribution is *usually* best. Pearson (1905, p. 181) introduced the excess when studying asymmetric laws.

In general, if the density is unknown, the knowledge of the first four moments is essential: when considering them as the corresponding theoretical moments of the density, it will be possible to imagine its type and therefore to calculate its parameters (hardly more than four of them).

To repeat: the normal law has only two parameters; therefore, if the calculated excess is sufficiently small, the unknown distribution will be determined by the first two moments. But what, indeed, is *sufficiently small*? We leave this question aside.

### 2.8. The Distribution of a Function of Random Variables

Suppose that random variables $\xi$ and $\eta$ have densities $\varphi_1(x)$ and $\varphi_2(y)$ and that $\eta = f(\xi)$ with a continuous and differentiable function $f$. The density $\varphi_1(x)$ is known and it is required to derive $\varphi_2(y)$. This problem is important and has to be solved often.

First of all, we (unnecessarily?) provide information about inverse functions and restrict our description to strictly monotone (increasing or decreasing) functions. The domain of an arbitrary function can however be separated into intervals of monotonic behaviour and each such interval can then be studied separately.

Suppose now that the function $y = f(x)$ strictly decreases on interval $[a; b]$. Turn its graph *to the left* until the $y$-axis is horizontal, and you will see the graph of the inverse function $x = \psi(y)$, also one-valued since the initial function was monotone. True, the positive direction of the $y$-axis and therefore of the argument $y$ (yes, $y$, not $x$ anymore) will be unusual. This nuisance disappears when looking with your mind's eye on the graph from the other side of the plane.

Return now to our problem. When $\xi$ moves along $[a; b]$, the random point $(\xi; \eta)$ moves along the curve $y = f(x)$. For example, if $\xi = x_0$, then $\eta = f(x_0) = y_0$. It is seen that the distribution function (not the density) $F(y)$ of $\eta$, or $P(\eta < y)$, is

$$P(\eta < y) = P(x < \xi < b) = \int_x^b \varphi_1(x) dx = \int_x^b \varphi_1(z) dz,$$

where $(x; y)$ is a current point on the curve $y = f(x)$.

Pursuing a methodical aim, we have changed the variable in the integral above but certainly did not alter the lower bound. However, it can be expressed as a function of $y$: $x = \psi(y)$. So now

$$F(y) = P(\eta < y) = \int_{\psi(y)}^b \varphi_1(z) dz.$$

Differentiate both parts of this equality with respect to $y$, and obtain thus the density

$$F'(y) = \varphi_2(y) = - \varphi_1[\psi(y)] \cdot \psi'(y).$$

For a strictly *increasing* function $f(x)$ the reasoning is the same although now it is the upper variable bound rather than the lower and the *minus* sign will disappear. Both cases can be written as

$$\varphi_2(y) = \varphi_1[\psi(y)] \cdot |\psi'(y)|.$$

*Example* (Ventzel 1969, p. 265).

$$\eta = 1 - \xi^3, \quad \varphi_1(x) = \frac{1}{\pi(1+x^2)}, \quad -\infty < x < \infty.$$

Here, $\varphi_1(x)$, is the Cauchy distribution (mentioned in § 2.3.2 in a slightly different form). We have

$$x = \psi(y) = \sqrt[3]{1-y},$$

$$\psi'(y) = -\frac{1}{3\sqrt[3]{(1-y)^2}}, \quad \varphi_1[\psi(y)] = f[\sqrt[3]{1-y}] = \frac{1}{\pi[1+\sqrt[3]{(1-y)^2}]},$$

$$\varphi_2(y) = \frac{1}{\pi[1+\sqrt[3]{(1-y)^2}]} \frac{1}{3\sqrt[3]{(1-y)^2}}.$$

Such a simple function … The $y$ can certainly be replaced by $x$.

### 2.9. The Bienaymé – Chebyshev Inequality

This is

$$P(|\xi - E\xi| < \beta) > 1 - \sigma^2/\beta^2, \quad \beta > 0 \qquad (2.24)$$

or, which is evident,

$$P(|\xi - E\xi| \geq \beta) < \sigma^2/\beta^2.$$

Inequality (2.24), and therefore its second form as well, take place for any random variable having an expectation and a variance and are therefore extremely interesting from a theoretical point of view. However, exactly this property means that the inequality is rather rough (I discuss any one of them). In a way, it combines the two magnitudes, $\sigma$ and $\beta$, without needing any other information.

Bienaymé (1853) established that inequality, but, unlike Chebyshev (1867 and later), did not pay special attention to his discovery since the subject of his memoir was not directly connected with it.

William Herschel (1817/1912, p. 579)
*presumed that any star promiscuously chosen […] out of* [more than 14 thousand] *is not likely to differ much from a certain mean size of them all.*

Stars unimaginably differ one from another and do not belong to a single population at all. The variance of their sizes is practically infinite, the notion of their mean size meaningless, and the inequality (2.24) cannot be applied. From another point of view, we may add: no positive data – no conclusion (*Ex nihilo nihil fit!*).

The English physician J. Y. Simpson (1847 – 1848/1871, p. 102) had similar thoughts: *The data* [about mortality after amputations] *have been objected to on the ground that they are collected from too many different hospitals and too many sources.* But […] *I believe […] that this very circumstance renders them more, instead of less, trustworthy.*

# Chapter 3. Systems of Random Variables. Correlation
## 3.1. Correlation

In the first approximation it may be stated that the variable *y* is a function of argument *x* on some interval or the entire number axis if, on that interval (on the entire axis), one and only one value of *y* corresponds to each value of *x*. Such dependence can exist between random variables. For example, Bessel (1838, §§ 1 – 2): the error of a certain type of measurements is $\eta = a\xi^2$.

Less tight connections between random variables are also possible (the stature of children depending on the stature of parents). Their study is the aim of an important chapter of mathematical statistics, of the theory of correlation. That word means *comparison*. More precisely, correlation considers the change *in the law of distribution* of a random variable depending on the change of another (or other) random variable(s) and as a rule on accompanying circumstances as well.

Lacking that specification and certainly without quantitative studies of phenomena (*qualitative*) correlation had been known in antiquity. (As stated in § 1.1.3, the entire ancient science had been qualitative.) Thus, Hippocrates (1952, No. 44): *Persons who are naturally very fat are apt to die earlier than those who are slender*. Climatic belts were isolated in antiquity, but only Humboldt (1817, p. 466) connected them with mean yearly air temperatures.

Seidel (1865 – 1866), a German astronomer and mathematician, first quantitatively investigated correlation. He studied the dependence of the monthly cases of typhoid fever on the level of subsoil water, and then both on that level and the rainfall.

Galton (1889) had begun to develop the theory of correlation proper, and somewhat later Pearson followed suit. Nevertheless, it had been sufficiently improved much later. Markov (1916/1951, p. 533) disparagingly but not altogether justly declared that the correlation theory's

*positive side is not significant enough and consists in a simple usage of the method of least squares for discovering linear dependences. However, not being satisfied with approximately determining various coefficients, the theory also indicates their probable errors, and enters here the realm of imagination* […].

Discovering dependences, even if only linear, is important and estimation of precision is certainly necessary. Linnik (Markov 1951, p. 670) noted that in those times correlation theory had still being developed so that Markov's criticism made sense. However, Hald (1998, p. 677), without mentioning either Markov or Linnik, described Fisher's pertinent contribution of 1915 (which Markov certainly did not see) and thus refuted Linnik. Anyway, here is Slutsky's reasonable general comment (letter to Markov of 1912, see Sheynin 1999b, p. 132):

*The shortcomings of Pearson's exposition are temporary and of the same kind as the known shortcomings of mathematics in the 17$^{th}$ and 18$^{th}$ centuries.*

Now we shall discuss the correlation coefficient. Two random variables, $\xi$ are $\eta$, are given. Calculate the moment

$$\mu_{\xi\eta} = E[(\xi - E\xi)(\eta - E\eta)] = E(\xi\eta) - 2E\xi\, E\eta + E\xi\, E\eta = E(\xi\eta) - E\xi\, E\eta$$

and divide it by the standard deviations $\sigma_\xi$ and $\sigma_\eta$ to obtain a dimensionless measure, the correlation coefficient

$$r_{\xi\eta} = \frac{\mu_{\xi\eta}}{\sigma_\xi \sigma_\eta}.$$

For independent $\xi$ and $\eta$ both $\mu_{\xi\eta}$ and (therefore) $r_{\xi\eta}$ disappear. The inverse statement is not true! Cf.: a sparrow is a bird, but a bird is not always a sparrow. One case is sufficient for refuting the inverse statement here also. So suppose that the density of $\xi$ is an even function, then $E\xi = 0$ and $E\xi^3 = 0$. Introduce now $\eta = \xi^2$, then $\mu_{\xi\eta} = E\xi^3 - 0 = 0$, QED. It follows that (even a functional) dependence can exist between random variables when the correlation coefficient is zero.

That coefficient takes values from $-1$ until 1. Correlation can therefore be negative. Example: the correlation (the dependence) between the stature and the weight of a person is positive, but between the distance from a lamp and its brightness is negative. Accordingly, we say that the correlation is direct or inverse.

### 3.2. The Distribution of Systems of Random Variables

Consider the probability $P(\xi < x, \eta < y)$. Geometrically, these inequalities correspond to an infinite region $-\infty < \xi < x, -\infty < \eta < y$, whereas analytically $P$ is expressed by the distribution function

$$F(x; y) = P(\xi < x, \eta < y).$$

Taken separately, the random variables $\xi$ and $\eta$ have distribution functions $F_1(x)$ and $F_2(y)$ and densities $f_1(x)$ and $f_2(y)$. Thus, for an infinitely large $x$ the inequality $\xi < x$ is identical and

$$F(+\infty; y) = F_2(y).$$

The density is also introduced here similar to the one-dimensional case:

$$P[(\xi; \eta) \text{ belongs to region } D] = \iint_D f(x; y)\,dxdy.$$

Function $f(x; y)$ is indeed the density. For independent $\xi$ and $\eta$ we have

$$f(x; y) = f_1(x) f_2(y).$$

For the bivariate (two-dimensional) normal law the density $f(x; y)$ is

$$\frac{1}{2\pi\sigma_x \sigma_y} \exp\{-\frac{1}{2(1-r^2)}[\frac{(x-E\xi)^2}{\sigma_\xi^2} - \frac{2r(x-E\xi)(y-E\eta)}{\sigma_\xi \sigma_\eta} + \frac{(y-E\eta)^2}{\sigma_\eta^2}]\}.$$

Apart from the previous notation, *r* is the correlation coefficient for ξ and η.

**3.2.1.** *The distribution of a sum of random variables.* Given, random variables ξ and η with densities $\varphi_1(x)$ and $\varphi_2(y)$. Required is the law of distribution of their sum $\omega = \xi + \eta$. For the distribution function of their system we have

$$F(x, y) = \iint \varphi(x, y) dx dy.$$

In case of infinite domains of both functions the integration is over an infinite half-plane

$$F(x; y) = \int_{-\infty}^{\infty} dx \int_{-\infty}^{\omega - x} \varphi(x; y) dy.$$

Differentiating it with respect to ω, we will have

$$F'_\omega(x; y) = \int_{-\infty}^{\infty} \varphi(x; \omega - x) dx = f(\omega),$$

or, after changing the places of *x* and *y*,

$$F'_\omega(x; y) = \int_{-\infty}^{\infty} \varphi(\omega - y; y) dy = f(\omega).$$

In case of independent random variables the distribution sought is called *composition* (of their densities). The formulas above lead to

$$f(\omega) = \int_{-\infty}^{\infty} \varphi_1(x) \varphi_2(\omega - x) dx = \int_{-\infty}^{\infty} \varphi_1(\omega - y) \varphi_2(y) dy.$$

Calculations are sometimes essentially simplified by geometrical considerations (Ventzel 1969, §§ 12.5 – 12.6). We only remark that the encounter problem (§ 1.1.2) can also be interpreted by means of the notions of random variable and density of distribution. Indeed, a random point (ξ; η) should be situated in a square with opposite vertices O(0, 0) and C(60, 60), and the sum (ξ + η) should be between two parallel lines, $y = x \pm 20$ (I have chosen 60 and 20 in that section). The distribution of that point could have ensured the derivation of the probability of the encounter. To recall, the moments ξ and η of the arrival of the two friends were independent.

## Chapter 4. Limit Theorems
### 4.1. The Laws of Large Numbers

The statistical probability $\hat{p}$ of the occurrence of an event can be determined by the results of independent trials, see formula (1.7), whereas its theoretical probability *p* is given by formula (1.1).

We (§ 1.1.1) listed the shortcomings of that latter formula and only repeat now that it is rarely applicable since equally probable cases are often lacking. Consequently, we have to turn to statistical probability.

**4.1.1.** *Jakob Bernoulli.* Bernoulli (1713/2005, pp. 29 – 30) reasonably remarked that

*Even the most stupid person […] feels sure that the more […] observations are taken, the less is the danger of straying from the goal.*

Nevertheless, he (p. 30) continued:

*It remains to investigate whether, when the number of observations increases, […] the probability of obtaining* [the theoretical probability] *continually augments so that it finally exceeds any given degree of certitude. Or* [to the contrary …] *that there exists such a degree of certainty which can never be exceeded no matter how the observations be multiplied.*

In other words, will the difference $|p - \hat{p}|$ continually decrease or not so that the statistical probability will not be a sufficiently good estimate of $p$.

Bernoulli proved that, in his restricted pattern, induction (trials) is (are) not worse than deduction: as $n \to \infty$ the difference $|p - \hat{p}|$ tends to disappear. His investigation opened up a new vast field. Nevertheless, not that difference itself, but its probability tends to zero. As $n \to \infty$

$$\lim P(|p - \hat{p}| < \varepsilon) = 1 \qquad (4.1)$$

with an arbitrarily small $\varepsilon$. This limit is exactly unity, not some lesser number, and induction is indeed not worse than deduction. But, wherever probabilities are involved, a fly appears in the ointment. The deviation of the statistical probability $\hat{p}$ from its theoretical counterpart can be considerable, even if rarely. A doubting Thomas can recall the example in § 1.2.3: the occurrence of an event with zero probability. Here, such an event is $|p - \hat{p}| \geq \varepsilon$. The limit of probability essentially differs from the *usual* limit applied in other branches of mathematics. Nothing similar can happen there!

Formula (4.1) is called the (weak) law of large numbers (a term due to Poisson). There also exists the so-called strong law of large numbers which removes the described pitfall, but I do not discuss it.

Strictly speaking, Bernoulli wrote out formula (4.1) in a somewhat different way, then continued his investigation. He proved that the inequality

$$P(|p - \hat{p}| \leq \varepsilon) \geq 1 - \delta, \delta > 0$$

will hold at given $\varepsilon$ and $\delta$ as soon as $n$ exceeds some $N$, which depends on those two numbers. He managed to determine how exactly $N$ must increase with the tightening of the initial conditions.

His investigation was not really successful: the demanded values of $N$ had later been considerably decreased (Pearson 1924; Markov 1924, p. 46 ff) mostly because it became possible to apply the Stirling formula unknown to Bernoulli. True, neither did De Moivre (§ 4.2) know that formula, but he derived it (even a bit before Stirling).

Mentioning Bernoulli's crude estimates Pearson (1925) inadmissibly compared his law with the wrong Ptolemaic system of the world. He missed its great general importance and, in particular, paid no attention to Bernoulli's existence theorem, of the very existence of the limit (4.1). It seems that Pearson did not set great store by such theorems.

In 1703 – 1705, before Bernoulli's posthumous *Ars Conjectandi* appeared, Bernoulli had exchanged letters with Leibniz; the Latin text of their correspondence are partially translated into German (Gini 1946; Kohli 1975); Bernoulli himself, without naming Leibniz, answered his criticisms in Chapter 4 of pt. 4 of his book. Leibniz did not believe that observations can ensure practical certainty and declared that the study of all the pertinent circumstances was more important than delicate calculations. Much later Mill (1843/1886, p. 353) supported this point of view:

*A very slight improvement of the data by better observations or by taking into fuller considerations the special circumstances of the case is of more use, than the most elaborate application of the calculus of probabilities founded on the* [previous] *data*.

He maintained that the neglect of that idea in applications to jurisprudence made the calculus of probability *the real opprobrium of mathematics*. Anyway, *considerations of the circumstances* and calculations do not exclude each other.

In a letter of 1714 to one of his correspondents Leibniz (Kohli 1975, p. 512) softened his doubts about the application of the statistical probability and mistakenly added that the late Bernoulli had *cultivated* [the theory of probability] in accordance with his, Leibniz', *exhortations*.

**4.1.2.** *Poisson*. Here is his qualitative definition of the law of large numbers (1837, p. 7):

*Les choses de toutes natures sont soumises à une loi universelle qu'on peut appeler <u>la loi des grands nombres</u>. Elle consiste en ce que, si l'on observe des nombres très considérables d'événements d'une même nature, dépendants de causes constantes et de causes qui varient irrégulièrement, tantôt dans un sens, tantôt dans l'autre, c'est-à-dire sans que leur variation soit progressive dans aucun sens déterminé, on trouvera, entre ces nombres, des rapports à très peu près constantes*.

This is a diffuse definition of a principle rather than law. And here is a contemporary qualitative definition of that law (Gnedenko 1954, § 30, p. 185): it is

*The entire totality of propositions stating with probability, arbitrarily close to unity, that there will occur some event depending on an indefinitely increasing number of random events each only slightly influencing it*.

The equality

$$\lim P(|\frac{\mu}{n} - \bar{p}| < \varepsilon) = 1, n \to \infty \qquad (4.2)$$

is now called the Poisson theorem. Here, µ/n is the frequency of an event in n independent trials and $p_k$ (from which $\bar{x}$ is calculated) is the probability of its occurrence in trial k. Note that for the *Bernoulli trials* the probability of the occurrence of the studied event was constant (not $p_k$ but simply p). Unlike formula (4.1), the new equality is general and therefore much more applicable.

**4.1.3.** *Subsequent history*. Chebyshev (1867) proved a more general theorem and Khinchin (1927) managed to generalize it still more. Finally, I provide another, not quite general formula for the law of large numbers: if

$$\lim P(|\bar{\xi}_n - a| < \varepsilon) = 1, n \to \infty,$$

where a is some number, the sequence of magnitudes $\xi_k$ obeys that law.

### 4.2. The De Moivre – Laplace Theorem

Suppose that a studied event occurs in each trial with probability p and does not occur with probability q, p + q = 1 and that in n such independent trials it happened µ times. Then, as $n \to \infty$,

$$\lim P(a \leq \frac{\mu - np}{\sqrt{npq}} \leq b) = \frac{1}{\sqrt{2\pi}} \int_a^b \exp(-\frac{z^2}{2}) dz. \qquad (4.3)$$

In the limit, the binomial distribution thus becomes normal. This is what De Moivre proved in 1733 for the particular case of p = q = 1/2 (in his notation, a = b = 1/2), but then he correctly stated that a transition to the general case is easy; furthermore, the heading of his (privately printed Latin) note mentioned the binomial $(a + b)^n$. To recall, np = Eµ and npq = varµ. Note also that the formula (4.3) is a particular case of the central limit theorem (§ 2.2.4).

Like other mathematicians of his time, De Moivre applied expansions into divergent series, only took into account their first terms, and neglected all the subsequent terms as soon as they began to increase (as soon as the series really began to diverge).

Laplace (1812, Chapter 3) derived the same formula (4.3) by means of a novelty, the Euler – Maclaurin summation formula. Furthermore, he added a term taking account of the inaccuracy occurring because of the unavoidable finiteness of n. Markov (1914/1951, p. 511), certainly somewhat mistakenly, called the *integral* after De Moivre and Laplace. That name persisted in Russian literature although, tacitly, in the correct way, as describing the integral theorem due to both those scholars. There also exists the corresponding local theorem

$$P(\mu) \approx \frac{1}{\sqrt{2\pi npq}} \exp[-\frac{(\mu - np)^2}{2npq}]. \qquad (4.4)$$

Assigning some µ in the right side of this formula and inserting the appropriate values of n, p, q, we will approximately calculate the probability of that µ. Exponential functions included in formulas (4.3) and (4.4) are tabulated in many textbooks.

A few additional remarks. Formula (4.3) describes a uniform convergence with respect to *a* and *b* (those interested can easily find this term), but Laplace (or certainly De Moivre) did not yet know that notion. Again, strict inequalities had not been then distinguished from non-strict ones. In formula (4.3), we should now apply a strict inequality in the second case (… < *b*). Then, the convergence to the normal law worsens with the decrease of *p* or *q* from 1/2 which is seen in a contemporary proof of the theorem (Gnedenko 1954, § 13).

In 1738 De Moivre included his own English translation of his private note in the second edition of the *Doctrine of Chances* and reprinted it in an extended form in the last edition (1756) of that book. However, the English language was not generally known by scientists on the Continent and the proof of (4.3) was difficult to understand since English mathematicians had followed Newton in avoiding the symbol of integral. Finally, Todhunter (1865, p. 192 – 193), the most eminent historian of probability of the 19th century, described the derivation of the formula (4.3) rather unsuccessfully and did not notice its importance. He even stated that De Moivre had only proved it for the particular case of *p* = *q* = 1/2. De Moivre's theorem only became generally known by the end of the 19th century.

Already in 1730 De Moivre independently derived the Stirling formula; the latter only provided him the value of the constant, $\sqrt{2\pi}$. Both Pearson (1924) and Markov (1924, p. 55 note) justly remarked that the Stirling formula ought to be called after both authors. I additionally remark that in 1730 De Moivre had compiled a table (with one misprint) of lg*n*! for *n* = 10(10)900 with 14 decimal points; 11 or 12 of them are correct.

Suppose now that it is required to determine the probability of casting a six 7 times in 100 rolls of a die. We have *p* = 1/6 and *n* = 100, then *np* = 16.7 and $\sqrt{npq}$ = 13.9. By formula (4.4)

$$P(\mu = 7) \approx \frac{1}{\sqrt{2\pi}\sqrt{npq}} \exp[-\frac{(7-np)^2}{2npq}].$$

I have isolated the factor $\sqrt{2\pi}$ since the exponential function is tabulated together with $1/\sqrt{2\pi}$.

Another point. In § 2.2.4 I mentioned that De Moivre had studied the sex ratio at birth. Now I say that exactly this subject (rather topical as the following shows) became the immediate cause for the derivation of formula (4.3).

Arbuthnot (1712) collected the data on births (or rather on baptisms) in London during 1629 – 1710. He noted that during each of those 82 years more boys had been born than girls and declared that that fact was *not the effect of chance, but Divine Providence, working for a good end* since mortality of boys and men was higher than that of females and since the probability of the observed fact was only $(1/2)^{-82}$.

His reasoning was not good enough but the problem itself proved extremely fruitful. Baptisms were not identical with births, London was perhaps an exception, Christians possibly somehow differed from other people and the comparative mortality of the sexes was not really studied. Then, by itself, an insignificant probability had not proved anything and it would have been much more natural to explain the data by a binomial distribution.

In a letter of 1713 Nikolaus Bernoulli (Montmort 1708/1713, pp. 280 – 285) had indeed introduced that distribution. Denote the yearly number of births by *n*, μ of them boys, the unknown sex ratio at birth by *m/f* and $p = m/(m + f)$. Bernoulli indirectly derived the approximate equality (lacking in Bernoulli's letter)

$$P(\frac{|\mu - np|}{\sqrt{npq}} \leq s) \approx 1 - \exp[-\frac{s^2}{2}],$$

where *s* had order √*n*, see Sheynin (1968; only in its reprint). He thus effectively arrived at the normal law much earlier than De Moivre.

Youshkevich (1986) reported that three mathematicians concluded that Bernoulli had come close to the local theorem (4.4) although I somewhat doubt it and the very fact that three mathematicians had to study Bernoulli's results testifies that these are difficult to interpret.

The initial goal of the theory of probability consisted in separating chance and design. Indeed, Arbuthnot, Nikolaus Bernoulli and De Moivre pursued this very aim. The last-mentioned devoted the first edition of his *Doctrine of Chances* to Newton and reprinted this dedication in the third edition of that book (p. 329). He attempted to work out, to *learn* from Newton's *philosophy*,

*A method of calculating the effects of chance* [… and to fix] *certain rules for estimating how far some sort of events may rather be owing to design rather than chance …*

Note that De Moivre then did not yet prove his limit theorem.

### 4.3. The Bayes Limit Theorem

His main formula was (1764)

$$P(b \leq r \leq c) = \int_b^c u^p (1-u)^q dx \div \int_0^1 u^p (1-p)^q dx. \qquad (4.5)$$

Bayes derived it by applying complicated logical constructions, but I interpret its conditions thus: given, a unit interval and segment [*b*; *c*] lying within it. Owing to complete ignorance (Scholium to Proposition 9), point *r* is situated with equal probability anywhere on that interval; in $n = p + q$ trials that point occurred *p* times within [*b*; *c*] and *q* times beyond it.

In other words, Bayes derived the posterior distribution of a random variable having a uniform prior distribution. The denial of that assumption (of the uniform distribution) led to discussions about the Bayes memoir (§ 1.1.1-5). In addition, the situation of point *r* is not at all random but unknown. Thus, the formula (4.5) should not be applied for deriving the probability of a certain value of a remote digit in the expansion of π (Neyman 1938a/1967, p. 337).

At that time there did not exist any clear notion of density; now, however, we may say that the formula (4.5) does not contradict its definition. Bayes derived the denominator of the formula and thus obtained the value of the beta-function (Euler). Both the pertinent calculation and the subsequent work were complicated and not easy to retrace. However, Timerding, the editor of the German version of Bayes' memoir (1908), surmounted the difficulties involved. Moreover, he invented a clever trick and wrote out the result as a limit theorem. For large *p* and *q* he arrived at

$$\lim P(a \leq \frac{\bar{p} - p/n}{\sqrt{pq/n^3}} \leq b) = \frac{1}{\sqrt{2\pi}} \int_a^b \exp(-\frac{z^2}{2}) \, dz, \, n \to \infty. \quad (4.6)$$

Here $\bar{p}$ is a statistical estimate of the unknown probability *p* that point *r* is within [*b*; *c*], and $p/n = \mathrm{E}\bar{p}$, $pq/n^3 = \mathrm{var}\,\bar{p}$.

A comparison of the formulas (4.3) and (4.5) convinces us that they describe the behaviour of differing random variables

$$\frac{\xi_i - \mathrm{E}\xi_i}{\sqrt{\mathrm{var}\xi_i}}, \, i = 1 \text{ (De Moivre)}, \, i = 2 \text{ (Bayes)}.$$

The variance in the Bayes formula is larger. The proof is not really needed; indeed, statistical data are present in both cases, but additional information (the theoretical probability) is only given in formula (4.3). And it is extremely interesting that Bayes, who had no idea about the notion of variance, understood that the De Moivre formula did not describe good enough the determination of that theoretical probability by its statistical counterpart. Both Jakob Bernoulli and De Moivre mistakenly stated the opposite, but Price, an eminent statistician who communicated (and extended) the posthumous Bayes memoir, mentioned this circumstance.

But why did not Bayes himself represent his result as a limit theorem? In another posthumous contribution of the same year, 1784, Bayes clearly indicated, for the first time ever, that the application of divergent series (in particular, by De Moivre) is fraught with danger. Timerding, it ought to be remarked, managed to avoid them. Note however that such series are still cautiously applied.

I believe that Bayes had completed the first version of the theory of probability which included the Jakob Bernoulli law of large numbers and the theorems due to De Moivre and Bayes himself. In addition, Bayes was actually the main predecessor of Mises (which the latter never acknowledged). See also Sheynin (2010b).

## Chapter 5. Random Processes. Markov Chains
### 5.1. Random Functions

Random functions are random variables changing discretely or continuously in time; for example, unavoidable noise occurring during the work of many instruments. Fixing some definite moment, we obtain the corresponding random variable, a *section of a random function*.

The law of distribution of a random function is naturally a function of two arguments one of which is time. For this reason the expectation of a random function is not a number but a (usual rather than a random) function. When fixing the moment of time, the expectation will pertain to the corresponding section of the random function and a similar statement concerns the variance. Another new point has to do with the addition of dependent random functions: the notion of correlation ought to be generalized.

A *random function without after-effect* is such for which there exists a definite probability of its transferring from a certain state to another one in such a way that additional information about previous situations does not change that probability. A good example is the Brownian motion (discovered by the English botanist Brown in 1827), the motion of tiny particles in a liquid under the influence of molecular forces.

About half a century ago, a new important phenomenon, the chaotic motion differing from random motion, began to be studied. However complicated and protracted is a coin toss, its outcomes do not change and neither do their probabilities. Chaotic motion, on the other hand, involves a rapid increase of its initial instability (of the unavoidable errors in the initial conditions) with time and countless positions of its possible paths.

It was Laplace (1781; 1812, § 15) who introduced subjective opinions (end of § 1.1.1-6) and, actually, a random process. Suppose that some interval is separated into equal or unequal parts and perpendiculars are erected from their ends. Let there be $i$ perpendiculars, their total length unity, forming a non-increasing sequence in one of the two directions. Their ends are connected by a broken line and a proper number of curves, and all this construction is repeated $n$ times after which the mean values of the current perpendiculars are calculated.

Laplace supposes that the lengths of the perpendiculars are assigned by $n$ different people and that the worth of candidates or the significance of various causes can thus be ordered in a non-increasing order. Each curve can be considered a random function; their set, a random process; and the mean curve, its expectation. True, the calculations occurred very complicated.

Evolution according to Darwin provides a second example. Consider a totality of individuals of, say, the male sex, of some species. Each individual can be theoretically characterised by the size of its body and body parts; the unimaginable multitude $n$ of such sizes is of no consequence. Introduce the usual definition of distance in an $n$-dimensional space and each individual will be represented by its point. The same space will contain the point or the subspace U of the sizes optimal for the chosen species. In the next generation, the offspring of any parents will be the better adapted to life the nearer they are to U which means that, in spite of the random scattering of the offspring around their midparents (a term due to Pearson), one generation after another will in general move towards U. However, that U will also move according to the changes in our surrounding world (and, if that movement is too rapid, the species can disappear). And so, individuals remote from U will in general perish or leave less offspring and our entire picture can be understood as a discrete random process with sections represented by each generation.

Our pattern is only qualitative; indeed, we do not know any numbers, any probabilities, for example, the probability of the mating of two given individuals of different sexes and we are therefore ignorant of any information about their offspring. Moreover, Darwin reasonably set great store by the habits of animals about which we are ignorant as well. Finally, there exists correlation between body parts of an individual. Darwin himself (1859/1958, p. 77) actually compared his theory (or, rather, hypothesis) with a random process:

*Throw up a handful of feathers, and all fall to the ground according to definite laws; but how simple is the problem where each shall fall compared with problems in the evolution of species.*

Opponents of evolution mostly cited the impossibility of its being produced by chance, by uniform randomness which once again showed that for a very long time that distribution had been considered as the only one describing randomness. Baer (1873, p. 6) and Danilevsky (1885, pt. 1, p. 194) independently mentioned the philosopher depicted in *Gulliver's Travels* (but borrowed by Swift form Raymond Lully, $13^{th} - 14^{th}$ centuries). That inventor, hoping to learn all the truths, was putting on record each sensible chain of words that appeared from among their uniformly random arrangements. Note that even such randomness does not exclude the gradual movement of the generations to U (but the time involved will perhaps be enormous).

Evolution began to be studied anew after Mendel's laws have been unearthed (after about 40 years of disregard) and, once more anew, after the important role of mutations has been understood.

### 5.2. Markov Chains (Processes with Discrete Time)

Suppose that one and only one of the events $A_1^{(s)}$, $A_2^{(s)}$, ..., $A_k^{(s)}$ occurs in trial *s* and that in the next trial the (conditional) probability of event $A_i^{(s+1)}$ depends on what happened in event *s*, but not on those preceding *s*. These conditions, if fulfilled in any trial, determine a *homogeneous Markov chain*.

Denote the conditional probability of $A_j^{(s+1)}$ as depending on $A_i^{(s)}$ by $p_{ij}$, then the process described by such chain is determined by a square matrix (table) of such probabilities, the *transition matrix*. Its first row is $p_{11}, p_{12}, \ldots, p_{1k}$, the second one, $p_{21}, p_{22}, \ldots, p_{2k}, \ldots$, and the last one, $p_{k1}, p_{k2}, \ldots, p_{kk}$, and the sum of the probabilities, *the transition probabilities*, in each row is unity.

It is possible to construct at once both a transition matrix for *n* trials and the limiting matrix which exists (that is, the corresponding limiting probabilities exist) if for some *s* all the elements of the matrix are positive. Markov derived this result and discovered some other findings which were later called ergodic theorems. In particular, it occurred that under certain conditions all the limiting probabilities are identical.

This remarkable property can explain, for example, the uniform distribution of the small planets along the ecliptic: a reference to these limiting probabilities which do not depend on the initial probabilities would have been sufficient. Actually, however, the small planets (more precisely, all planets) move along elliptical orbits and in somewhat differing planes.

Poincaré (1896/1987, p. 150), who had not referred to any Russian author, not even to Laplace or Poisson, justified this fact although in a complicated way. (Also, by introducing hypercomplex numbers, he proved that after a lot of shuffling the positions of the cards in a pack tended to become equally probable.)

Markov himself only applied his results to investigate the alternation of consonants and vowels in the Russian language (Petruszewycz 1983). He possibly obeyed his own restriction (Ondar 1977/1981, p. 59, Markov's letter to Chuprov of 1910):

*I shall not go a step out of that region where my competence is beyond any doubt*.

The term itself, *Markov chain*, first appeared (in French) in 1926 (Bernstein 1926, first line of § 16) and pertained to Markov's investigations of 1906 – 1913. Some related subjects are Brownian motion, extinction of families, financial speculation, random walk.

The urn problem discussed below can be understood as a (one-dimensional) random walk, as a discrete movement of a particle in one or another direction along some straight line with the probabilities $p$ and $q$ of movement depending on what had happened in the previous discrete moment. Diffusion is a similar but three-dimensional process, but a random walk with constant $p$ and $q$, like the *walk* of the number of winnings of one of the two gamblers in a series of games, is not anymore a Markov chain.

And so, we will discuss the urn problem of Daniel Bernoulli (1770) and Laplace which is identical to the celebrated Ehrenfests' model (1907) considered as the beginning of the history of discrete random processes, or Markov chains. The first urn contains $n$ white balls, the second urn, the same number of black balls. Required is the (expected) number of white balls in the first urn after $r$ cyclic interchanges of a ball.

In his second problem Bernoulli generalized the first by considering three urns and balls of three colours. He managed to solve it elegantly, and discovered the limiting situation, an equal (expected) number of balls of each colour in each urn. A simplest method of confirming this result consists in a reference to the appropriate ergodic theorem for homogeneous Markov chains, but first we should prove that this Bernoulli problem fits the pattern of that theorem. It is not difficult. Indeed, for example, in the case of two urns, four events are possible at each interchange and the probability of each event changes depending on the results of the previous interchange. These four events are: white balls were extracted from each urn; a white ball was extracted from the first urn and a black ball from the second etc.

Laplace (1812, chapter 3) generalized the Bernoulli problem (but did not refer to him) by admitting an arbitrary initial composition of the urns, then (1814/1886, p. LIV) adding that *new urns are placed*

*amongst the original urns*, again with an arbitrary distribution of the balls. He (p. LV) concluded, probably too optimistically, that

*On peut étendre ces résultats à toutes les combinaisons de la nature, dans lesquelles les forces constantes dont leurs éléments sont animés établissent des modes réguliers d'action, propres à faire éclore du sein même du chaos des systèmes régis des lois admirables.*

## Chapter 6. The Theory of Errors and the Method of Least Squares
### 6.1. The Theory of Errors

This term (in German) is due to Lambert (1765, § 321). It only became generally used in the middle of the next century; neither Laplace, nor Gauss ever applied it although Bessel did. The theory of errors studies errors of observation and their treatment so that the method of least squares (MLSq) belongs to it. I have separated that method owing to its importance.

The theory of errors can be separated into a stochastic and a determinate part. The former studies random errors and their influence on the results of measurements, the action of the round-off errors and, the dependence between obtained relations. The latter investigates the patterns of measurement for a given order of errors and studies methods of excluding systematic errors (or minimizing their influence).

Denote a random error by $\xi$, its expectation will then be $E\xi = 0$. Otherwise (as it really is) $\xi$ is not a purely random error and $E\xi = a$ is the systematic error. It shifts the even density of random errors either to the right (if $a > 0$), or to the left (if $a < 0$).

From 1756 (Simpson, § 1.2.3) until the 1920s the stochastic theory of errors, as stated there, had remained the most important field of application of the theory of probability. In a posthumous publication Poincaré(1921/1983, p. 343) noted that *La théorie des erreurs était naturellement mon principal but* and Lévy (1925, p. vii) strongly indicated that without the theory of errors his contribution on stable laws *n'aurait pas de raison d'être*.

Stable laws became an essential notion of the theory of probability, but for the theory of errors they are absolutely useless. As a corollary to the definition of a stable law it follows that the sum $\sum \xi_i$ and the arithmetic mean $\bar{\xi}$ have the same distribution as the independent and identically distributed random variables $\xi_i$, and Lévy proved that a real estimation of the precision of those functions of random variables, if their distribution is not stable, is very difficult. However, an observer can never know whether the errors of his measurements obey a stable law or not. Moreover, the Cauchy law is also stable, but does not possess any variance (§ 2.3.2).

In turn, mathematical statistics took over the principles of maximal likelihood and least variance (see § 6.3 below) from the stochastic theory of errors.

Now the determinate theory of errors. Hipparchus and Ptolemy could not have failed to be ignorant about them (in the first place, about those caused by the vertical refraction). Nevertheless, it was Daniel Bernoulli (1780) who first clearly distinguished random and systematic errors although only in a particular case.

Also in antiquity astronomers had been very successfully observing under the most favourable conditions. A good example is an observation of the planets during their *stations*, that is, during the change of their apparent direction of motion, when an error in registering some definite moment least influences the results of subsequent calculations. Indeed,

*One the most admirable features of ancient astronomy* [was] *that all efforts were concentrated upon reducing to a minimum the influence of the inaccuracy of individual observations with crude instruments by developing* […] *the mathematical consequences of very few elements* [of optimal circumstances] (Neugebauer 1950/1983, p. 250).

Actually, however, the determinate theory of errors originated with the differential calculus. Here is a simplest geodetic problem. Two angles, α and β, and side *a* are measured in a triangle and the order of error of these elements is known. Required is the order of error in the other (calculated) elements of the triangle, and thus the determination of the optimal form of the triangle.

Denote the length of any of the calculated sides by *W*. It is a function of the measurements:

$W = f(a; α; β)$,

and its differential, approximately equal to its error, is calculated by standard formulas.

From studying isolated geodetic figures the determinate theory moved to investigating chains and even nets of triangles. And here is a special problem showing the wide scope of that theory (Bessel 1839). A measuring bar several feet in length is supported at two points situated at equal distances from its middle. The bar's weight changes its length and the amount of change depends on the position of the supporting points. Where exactly should you place these points so that the bar's length will be least corrupted? Bessel solved this problem by means of appropriate differential equations. For a contemporary civil engineer such problems are usual, but Bessel was likely the first in this area.

Gauss and Bessel originated a new stage in experimental science. Indeed, Newcomb (Schmeidler 1984, pp. 32 – 33) mentioned the *German school of practical astronomy* but mistakenly only connected it with Bessel. True, the appropriate merits of Tycho Brahe are not known adequately. Newcomb continued:

*Its fundamental idea was that the instrument is indicted* […] *for every possible fault, and not exonerated till it has proved itself corrected in every point. The methods of determining the possible errors of an instrument were developed by Bessel with ingenuity and precision of geometric method …*

Gauss had detected the main systematic errors of geodetic measurements (those caused by lateral refraction, by the errors of graduating the limbs of the theodolites, and inherent in some methods of measurement) and outlined the means for eliminating/decreasing them.

For a more detailed description of this subject see Sheynin (1996, Chapter 9).

## 6.2. The True Value of a Measured Constant

Many sciences and scientific disciplines have to measure constants; metrology ought to be mentioned here in the first place. But what should we understand as a true value of a constant? Is it perhaps a philosophical term?

Fourier (1826) suggested its definition undoubtedly recognized earlier even if tacitly: the true value of a constant is the limit of the arithmetic mean of its $n$ measurements as $n \to \infty$. It is easy to see that the Mises' frequentist definition of probability (§ 1.1.3) is akin to Fourier's proposal. The measurements can be carried out under identical conditions as is really essential for metrology but this is inadmissible in geodesy: differing (but good enough) external conditions were necessary for some compensation of systematic errors.

I failed to find a single reference to his definition but many authors repeated it independently from him or one another. One of them (Eisenhart 1963/1969, p. 31) formulated the unavoidable corollary: the mean residual systematic error had to be included in the *true value*:

*The mass of a mass standard is […] specified […] to be the mass of the metallic substance of the standard plus the mass of the average volume of air adsorbed upon its surface under standard conditions.*

Statisticians have done away with true values and introduced instead parameters of densities (or distribution functions) and their properties. A transition to more abstract notions is a step in the right direction (cf. end of § 1.2.3), but they still have to mention true values; Gauss (1816, §§ 3 и 4) even discussed the true value of a measure of error, of something not existing in nature. For more detail see Sheynin (2007).

## 6.3. The Method of Least Squares

This standard method of treating observations is usually regarded as a chapter of mathematical statistics rather than probability.

Suppose that the unknown constants $x$ and $y$ are connected with observations $w_1, w_2, \ldots, w_n$ by linear equations

$$a_i x + b_i y + \ldots + w_i = 0, \quad i = 1, 2, \ldots, n. \qquad (6.1)$$

In the general case the number of the unknowns, $k$, is arbitrary, but if $k > n$, the solution of (6.1) is impossible, and if $k = n$, no special methods of its solution are needed. Therefore, $k < n$. The coefficients $a_i, b_i, \ldots$ ought to be provided by the appropriate theory, and the system (6.1) can be supposed linear if the unknowns are small.

Indeed. Suppose that a system is not linear and that its first equation is

$$a_1 x^2 + b_1 y^2 + w_1 = 0.$$

We actually know the approximate value of the unknowns, $x_0$ and $y_0$, so that, for example, the first term of our equation is $a_1(x_0 + \Delta x)^2$ with an unknown small $\Delta x$. The term $(\Delta x)^2$ can be neglected and that first term will be $a_1(x_0^2 + 2x_0\Delta x)$. The unknown magnitude is now linear, $2a_1 x_0 \Delta x$, and the second term of our equation, $b_1 y^2$, can be *linearized* in a similar way. A similar procedure is easily explained in case of trigonometric coefficients which is important for astronomy.

And now the main question: how to solve the system (6.1)? Observations are supposed to be independent (or almost so) and rejecting the $(n - k)$ redundant equations (which exactly?) would have been tantamount to rejecting worthy observations. A strict solution is impossible: any set $(x, y, \ldots)$ will leave some residual free terms (call them $v_i$). We are therefore obliged to impose one or another condition on these residuals. It became most usual to choose the condition of least squares

$$v_1^2 + v_2^2 + \ldots + v_n^2 = \min, \qquad (6.2)$$

hence the name, MLSq. And

$$v_i^2 = (a_i x + b_i y + \ldots + w_i)^2.$$

We ought to determine the values of $x, y, \ldots$, leading to condition (6.2), and these unknowns are therefore considered here as variables. We have

$$\frac{\partial v_i^2}{\partial x} = 2a_i(a_i x + b_i y + \ldots + w_i).$$

According to the standard procedure,

$$\frac{\partial v_1^2}{\partial x} + \frac{\partial v_2^2}{\partial x} + \ldots + \frac{\partial v_n^2}{\partial x} = 2\sum(a_i a_i x + a_i b_i y + \ldots + a_i w_i) = 0,$$

so that, applying the Gauss elegant notation (§ 2.3.2),

$$[aa]x + [ab]y + \ldots + [aw] = 0.$$

Differentiating $v_i^2$ with respect to $y$, we similarly get

$$[ab]x + [bb]y + \ldots + [bw] = 0.$$

The derived equations are called *normal*, and it is clear that their number coincides with the number of the unknowns (yes, they became again unknown); the system of normal equations can therefore be solved in any usual way. Note, however, that the solution provides a certain set $\hat{x}, \hat{y},\ldots$, a set of estimators of $\{x, y, \ldots\}$, of magnitudes which will remain unknown. Even the unknowns of the system of normal equations already are $\hat{x}, \hat{y},\ldots$ rather than $x, y, \ldots$ It is also necessary to estimate the errors of $\hat{x}, \hat{y},\ldots$, but we leave that problem aside.

Condition (6.2) ensures valuable properties to those estimators (Petrov 1954). It corresponds to the condition of minimal variance, to

$$m^2 = \frac{v_1^2 + v_2^2 + ... + v_n^2}{n-k} = \min. \tag{6.3}$$

The denominator is the number of redundant observations; the same is true for the formula (2.12) which corresponds to the case of one single unknown. For this case system (6.1) becomes simpler,

$a_i x + w_i = 0,$

and it is easy to verify that it leads to the generalized arithmetic mean.

Classical systems (6.1) had two unknown parameters of the ellipsoid of rotation best representing the figure of the Earth. After determining the length of one degree of a meridian in two different and observed latitudes it became possible to calculate those parameters whereas redundant *meridian arc measurements* led to equations of the type of (6.1). They served as a check of field measurements, they also heightened the precision of the final result (and to some extent compensated local irregularities of the figure of the Earth).

The lengths of such arcs in differing latitudes were certainly different and thus indicated the deviation of that figure from a circumference.

Legendre (1805, pp. 72 – 73; 1814) recommended the MLSq although only justifying it by reasonable considerations. Moreover, he (as also Laplace) mistakenly called the $v_i$'s errors of measurements and, finally, according to the context of his recommendation, he thought that the MLSq led to the minimal value of the maximal $|v_i|$. Actually, this condition is ensured by the method of *minimax*, see § 6.4.

Gauss had applied the MLSq from 1794 or 1795 and mistakenly thought that it had been known long ago. In general, Gauss did not hurry to publish his discoveries; he rather connected priority with the finding itself. He (1809, § 186) therefore called the MLSq *unser Princip* which deeply insulted Legendre. Note, however, that, unlike Legendre, Gauss had justified the new method (but later substantiated it in a different way since then, in 1809, the normal law became the only law of error).

As I see it, Legendre could have simply stated in some subsequent contribution that no one will agree that Gauss was the author of the MLSq. However, French mathematicians including Poisson (see below a few words about Laplace) sided with Legendre's opinion and, to their own great disadvantage, ignored Gauss' contributions on least squares and the theory of errors.

Much later Gauss (letter to Bessel of 1839; *Werke*, Bd. 8, pp. 146 – 147) explained his new attitude towards the MLSq:

*Ich muß es nämlich in alle Wege für weniger wichtig halten, denjenigen Wert einer unbekannten Größe auszumitteln, dessen Wahrscheinlichkeit die größte ist, die ja doch immer unendlich klein bleibt, als vielmehr denjenigen, an welchen sich haltend man das am wenigsten nachteilige Spiel hat.*

In other words, an integral measure of reliability (the variance) is preferable to the principle of maximal likelihood which he applied in 1809. Then, in 1809, Gauss did not refer either to Lambert (1760, §

295) or to Daniel Bernoulli (1778). The former was the first to recommend that principle for an indefinite density distribution. He had only graphically shown that density; it was a more or less symmetrical unimodal curve of the type φ($x - \hat{x}$), where $\hat{x}$ can be understood as a location parameter. For observations $x_1, x_2, …, x_n$ Lambert recommended to derive $\hat{x}$ from the condition (of maximum likelihood nowadays applied in mathematical statistics)

φ($x_1 - \hat{x}$) φ($x_2 - \hat{x}$) … φ($x_n - \hat{x}$) = max.

So Gauss assumed that the arithmetic mean of observations was at the same time the most probable (in the sense of maximum likelihood) estimator and discovered that only the normal distribution followed from this assumption.

In 1823 Gauss published his second and final justification of the MLSq by the principle of minimal variance (see above his letter to Bessel of 1839). Unlike his considerations in 1809, his reasoning which led him to equations (6.2) was very complicated whereas the law of error was indeed more or less normal. Thus, Maxwell (1860) proved (non-rigorously) that the distribution of gas velocities appropriate to a gas in equilibrium was normal; Quetelet (1853, pp. 64 – 65) maintained that the normal law governed the errors *faites par la nature*.

No wonder that Gauss' first formulation of the MLSq persisted (and perhaps is still persisting) in spite of his opinion. I (2012) noticed, however, that Gauss actually derived formula (6.3) as representing the minimal value of the (sample) variance independently from his complicated considerations and that, when taking this into account, his memoir (1823) becomes much easier to understand. And facts showing that the normal law was not universal in nature continued to multiply so that that memoir should be considered much more important.

The first serious opponent of the normal law was Newcomb (1886, p. 343) who argued that the cases of normal errors were quite exceptional. For treating long series of observations he recommended a mixture of differing normal laws, but the calculations proved complicated whereas his pattern involved subjective decisions. Later Eddington (1933, § 5) proved that that mixture was not stable.

Bessel (1818) discussed Bradley's series of 300 observations and could have well doubted the existence of an appropriate normal law. He noticed that large errors had appeared there somewhat more often than expected but somehow explained it away by the insufficient (!) number of observations. Much later he (1838) repeated his mistake. I (2000) noted many other mistakes and even absurdities in his contributions.

Unlike other French mathematicians, Laplace objectively described Legendre's complaint: he was the first to publish the MLSq, but Gauss had applied it much earlier. However, Laplace never recognized the utmost importance of Gauss' second substantiation of the method. Instead, he persisted in applying and advocating his own version of substantiating it. He proved several versions of the central limit theorem (CLT) (§ 2.2.4), certainly, non-rigorously (which was quite understandable) and very carelessly listing its conditions, then declared that the errors of observation were therefore normal. Laplace (1814/1886, p. LX) maintained that his finding was applicable in astronomy where long series of observations are made; cf., however, Newcomb's opinion above. Then he (1816/1886, p. 536) stated that the CLT holds in geodesy since, as it followed from his reasoning, the order of two main errors inherent in geodetic observations have been equalized. Here again he did not really take into account the conditions of that theorem.

Markov's work on the MLSq has been wrongly discussed. Neyman (1934) attributed to him Gauss' second justification of 1823 which even until our time (Dodge 2003, p. 161) is sometimes called after both Gauss and Markov. David & Neyman (1938) repeated the latter's mistake but the same year Neyman (1938b/1952, p. 228) corrected himself.

Then, Linnik et al (1951, p. 637) maintained that Markov had *in essence* introduced concepts *identical* to the modern notions of unbiased and effective statistics. Without explaining that latter notion I simply note that these authors should have replaced Markov by Gauss.

Markov (1899) upheld the second justification of the MLSq perhaps much more resolutely than his predecessors (the first such opinion appeared in 1825). However, he (1899/1951, p. 246) depreciated himself:

*I consider* [that justification] *rational since it does not obscure the conjectural essence of the method.* […] *We do not ascribe the ability of providing the most probable or the most plausible results to the method …*

Such a method does not need any justification. Furthermore, the MLSq does have optimal properties (Petrov 1954, cited above as well). Also see Sheynin (2006).

### 6.4. The Minimax Method

There also exist other methods of solving systems (6.1). They do not lead to the useful properties of the MLSq estimators but are expedient in some cases. And there also exists a special method leading to the least absolute value of the maximal residual

$$|v_{max}| = \min. \tag{6.4}$$

Least means least among all possible solutions (and therefore sets of $v_i$'s); in the simplest case, among several reasonable solutions. The minimax method does not belong to probability theory, does not lead to any *best* results, but it allows to make definite conclusions. Recall that the coefficients $a_i$, $b_i$, … in system (6.1) are given by the appropriate theory and ask yourselves: do the observations $w_i$ confirm it? After determining $\hat{x}, \hat{y},...$ (this notation does not infer the MLSq anymore) we may calculate the residual free terms $v_i$ and determine whether they are not too large as compared with the known order of errors. In such cases we ought to decide whether the theory was wrong or that the observations were substandard. And here the minimax method is important: if even condition (6.4) leads to inadmissible $|v_i|$, our doubts are certainly justified.

Both Euler and Laplace had applied the minimax method (the latter had devised an appropriate algorithm) for establishing whether the accomplished meridian arc measurements denied the ellipticity of the figure of the Earth. Kepler (1609/1992, p. 334/143) could have well applied elements of that method for justifying his rejection of the Ptolemaic system of the world: the Tychonian observations were sufficiently precise but did not agree with it. In astronomy, equations are neither linear, nor even algebraic, and Kepler had to surmount additional difficulties (irrespective of the method of their solution).

Condition (6.4) is identical to having

$$\lim(v_1^{2k} + v_2^{2k} + ... + v_n^{2k}) = \min, k \to \infty,$$

which is almost obvious. Indeed, suppose that $|v_1|$ is maximal. Then, as $k \to \infty$, all the other terms of the sum can be neglected. For arriving at a minimal value of the sum, $v_1^{2k}$, and therefore $|v_1|$ also, should be as small as possible.

Without looking before he leapt, Stigler(1986, pp. 27, 28) confidently declared that Euler's work (see above) was *a statistical failure* since he

*Distrusted the combination of equations, taking the mathematician's view that errors actually increase with aggregation rather than taking the statistician's view that random errors tend to cancel one another.*

However, at the turn of the 18<sup>th</sup> century Laplace, Legendre and other scientists *refusa de compenser* les angles of a triangulation chain between two bases. Being afraid of corrupting the measured angles by adjusting them, they resolved to calculate each half of the chain from its own base (and somehow to adjust the common side of both parts of the triangulation), see Méchain & Delambre (1810, pp. 415 – 433). Later, in the Third Supplement to his *Théorie analytique*, Laplace (ca. 1819/1886, pp. 590 – 591) explained that decision by the lack of the *vraie théorie* which he (rather than Gauss whom he had not mentioned) had since created. See also Sheynin (1993a, p. 50).

In the Soviet Union, separate triangulation chains were included in the general adjustment of the entire network only after preliminarily adjustment (§ 1.1.4). This pattern was necessary since otherwise the work would have been impossible. In addition, the influence of the systematic errors should have been restricted to separate chains (as stated in a lecture of ca. 1950 of an eminent Soviet geodesist, A. A. Isotov), and this consideration was akin to the decision of the French scientists described above.

### Chapter 7. Theory of Probability, Statistics, Theory of Errors
### 7.1. Axiomatization of the Theory of Probability

Following many previous author, I noted (§ 1.1.1) that the classical definition of probability is unsatisfactory and that Hilbert (1901/1970, p. 306) recommended to axiomatize the theory of probability:

*Durch die Untersuchungen über die Grundlagen der Geometrie wird uns die Aufgabe nahe gelegt, nach diesem Vorbilds diejenigen physikalischen Disziplinen axiomatisch zu behandeln, in denen schon*

*heute die Mathematik eine hervorragende Rolle spielt: dies sind in erster Linie die Wahrscheinlichkeitsrechnung und die Mechanik.*

The theory of probability had then still been an applied mathematical (but not physical) discipline. In the next lines of his report Hilbert mentioned the method of mean values. That method or theory had been an intermediate entity divided between statistics and the theory of errors, and Hilbert was one of the last scholars (the last one?) to mention it, see Sheynin (2007, pp. 44 – 46).

Boole (1854/1952, p. 288) indirectly forestalled Hilbert:

*The claim to rank among the pure sciences must rest upon the degree in which it* [the theory of probability] *satisfies the following conditions: 1° That the principles upon which its methods are founded should be of an axiomatic nature.*

Boole formulated two more conditions, both of a general scientific nature. Attempts to axiomatize the theory began almost at once after Hilbert's report. However, as generally recognized, only Kolmogorov attained quite satisfactory results. Without discussing the essence of his work (see for example Gnedenko 1954, § 8 in chapter 1), I quote his general statements (1933, pp. III and 1):

*Der leitende Gedanke des Verfassers war dabei, die Grundbegriffe der Wahrscheinlichkeitsrechnung, welche noch unlängst für ganz eigenartig galten, natürlicherweise in die Reihe der allgemeinen Begriffsbildungen der modernen Mathematik einzuordnen.*

*Die Wahrscheinlichkeitstheorie als mathematische Disziplin soll und kann genau in demselben Sinne axiomatisiert werden wie die Geometrie oder die Algebra. Das bedeutet, daß, nachdem die Namen der zu untersuchenden Gegenstände und ihrer Grundbeziehungen sowie die Axiome, denen diese Grundbeziehungen zu gehorchen haben, angegeben sind, die ganze weitere Darstellung sich ausschließlich auf diese Axiome gründen soll und keine Rücksicht auf die jeweilige konkrete Bedeutung dieser Gegenstände und Beziehungen nehmen darf.*

For a long time these ideas had not been generally recognized (Doob 1989; 1994, p. 593):

*To most mathematicians mathematical probability was to mathematics as black marketing to marketing;* […] *The confusion between probability and the phenomena to which it is applied* […] *still plagues the subject;* [the significance of the Kolmogorov monograph] *was not appreciated for years, and some mathematicians sneered that* […] *perhaps probability needed rigor, but surely not rigor mortis;* […] *The role of measure theory in probability* […] *still embarrasses some who like to think that mathematical probability is not a part of analysis.*

*It was some time before Kolmogorov's basis was accepted by probabilists. The idea that a (mathematical) random variable is simply a function, with no romantic connotation, seemed rather humiliating to some probabilists …*

For a long time Hausdorff's merits had remained barely known. His treatise on the set theory (1914, pp. 416 – 417) included references to probability, but much more was contained in his manuscripts, see

Girlich (1996) and Hausdorff (2006). I also mention Markov (1924). On p. 10 he stated a curious axiom and on p. 24 referred to it (without really thinking how the readers will manage to find it):

*Axiom.* [Not separated from general text!] *If* [...] *events p, q, r, ..., u, v are equally possible and divided with respect to event A into favourable and unfavourable, then*, [if] *A has occurred,* [those] *which are unfavourable to event A fall through, whereas the others remain equally possible ...*

*The addition and multiplication theorems along with the axiom mentioned above serve as an unshakeable base for the calculus of probability as a chapter of pure mathematics*.

His axiom and statement have been happily forgotten although he had a predecessor. Donkin (1851) had introduced a similar axiom but did not claim to change the status of probability theory and Boole (1854/2003, p. 163) positively mentioned him.

Shafer & Vovk (2001) offered their own axiomatization, possibly very interesting but demanding some financial knowledge. They (2003, p. 27) had explained their aim:

[In our book] *we show how the classical core of probability theory can be based directly on game-theoretic martingales, with no appeal to measure theory. Probability again becomes* [a] *secondary concept but is now defined in terms of martingales ...*

Barone & Novikoff (1978) and Hochkirchen (1999) described the history of the axiomatization of probability. The latter highly estimated an unpublished lecture of Hausdorff read in 1923.

### 7.2. Definitions and Explanations

As understood nowadays, statistics originated in political arithmetic (Petty, Graunt, mid-17$^{th}$ century). It quantitatively (rather than qualitatively) studied population, economics and trade, discussed the appropriate causes and connections and applied simplest stochastic considerations. Here is a confirmation (Kendall 1960):

*Statistics, as we now understand the term, did not commence until the 17$^{th}$ century, and then not in the field of 'statistics'* [Staatswissenschaft]. *The true ancestor of modern statistics is* [...] *Political Arithmetic*.

Statistics had gradually, and especially since the second half of the 19$^{th}$ century, begun to penetrate various branches of natural sciences. This led to the appearance of the term *statistical method* although we prefer to isolate three stages of its development.

At first, conclusions were being based on (statistically) noticed qualitative regularities, a practice which conformed to the qualitative essence of ancient science. See the statements of Hippocrates (§ 3.1) and Celsus (§ 1.1.3).

The second stage (Tycho in astronomy, Graunt in demography and medical statistics) was distinguished by the availability of statistical data. Scientists had then been arriving at important conclusions either by means of simple stochastic ideas and methods or even directly, as before. A remarkable example is the finding of an English physician Snow (1855/1965, pp. 58 – 59) who compared mortality from cholera for two groups of the London population, of those whose drinking water was (somehow) purified or not. Purification decreased mortality by 8 times and he thus discovered the way in which cholera epidemics had been spreading.

During the present stage, which dates back to the end of the 19th century, inferences are being checked by quantitative stochastic rules.

The questions listed by Moses (Numbers 13:17 – 20) can also be attributed to that first stage (and to political arithmetic): he sent scouts to *spy out* the land of Canaan, to find out

*whether the people who dwell in it are strong or weak, whether they are few or many*, […] *whether the land is rich or poor …*

In statistics itself, exploratory data analysis was isolated. Already Quetelet discussed its elements (1846); actually, however, the introduction of isolines was a most interesting example of such analysis. Humboldt (1817, p. 466) invented isotherms and (much later) mentioned Halley who, in 1791, had shown isolines of magnetic declination over North Atlantic.

That analysis belongs to the scientific method at large rather than mathematics and is not therefore recognized by mathematical statistics. It only belongs to theoretical statistics which in my opinion should mostly explain the difference between the two statistical sisters. Some authors only recognize either one or another of them. In 1953 Kolmogorov (Anonymous 1954, p. 47), for example, declared that

*We have for a long time cultivated a wrong belief in the existence, in addition to mathematical statistics and statistics as a socio-economic science, of something like yet another non-mathematical, although universal <u>general</u> theory of statistics which essentially comes to mathematical statistics and some technical methods of collecting and treating statistical data. Accordingly, mathematical statistics was declared a part of this <u>general theory of statistics</u>. Such views […] are wrong.* […]

*All that which is common in the statistical methodology of the natural and social sciences, all that which is here indifferent to the specific character of natural or social phenomena, belongs to* […] *mathematical statistics.*

These *technical methods* indeed constitute exploratory data analysis.

Kolmogorov & Prokhorov (1974/1977, p. 721) defined mathematical statistics as

*the branch of mathematics devoted to the mathematical methods for the systematization, analysis and use of statistical data for the drawing of scientific and practical inferences.*

Recall (§ 2.7) that they also defined the notion of statistical data and note that a similar definition of the theory of statistics had appeared in the beginning of the 19th century (Butte 1808, p. XI): it is

*Die Wissenschaft der Kunst statistische Data zu erkennen und zu würdigen, solche zu sammeln und zu ordnen.*

The term *mathematical statistics* appeared in the mid-19th century (Knies 1850, p. 163; Vernadsky 1852, p. 237), and even before Butte Schlözer (1804) mentioned the theory of statistics in the title of his book. He (p. 86) also illustrated the term *statistics*: *Geschichte ist eine fortlaufende Statistik, und Statistik stillstehende Geschichte.*

Obodovsky (1839, p. 48) offered a similar statement: *history is related to statistics as poetry to painting*.

Unlike Shlözer, many statisticians understood his pithy saying as a *definition* of statistics; as well we may say today: a car is a landed plane, and a plane, a car taken wing. Schlözer had not noticed that statistics ought to be compiled in a certain state or locality at differing periods of time and compared with each other. Indeed, this is extremely important so that statistics is not standing still!

For us, the theory of statistics essentially originated with Fisher. A queer episode is connected here with Chuprov's book (1909/1959). Its title is *Essays on the Theory of Statistics*, but on p. 20 he stated that *A clear and strict justification of the statistical science is still needed*!

The determinate theory of errors (§ 6.1) has much in common with both the exploratory data analysis and Fisher's creation, the experimental design (a rational organization of measurements corrupted by random errors). However, the entire theory of errors seems to be the application of the statistical method to the process of measurement and observation in experimental science rather than a chapter of mathematical statistics, as it is usually maintained. Indeed, stellar statistics is the application of the statistical method to astronomy, and medical statistics is etc. Furthermore, unlike mathematical statistics the theory of errors cannot *at all* give up the notion of true value of a measured constant (§ 6.2).

### 7.3. Penetration of the Theory of Probability into Statistics

Hardly anyone will deny nowadays that statistics is based on the theory of probability, but the situation had not always been the same. Already Jakob Bernoulli (§ 4.1.1) firmly justified the possibility of applying the latter but statisticians had not at all been quick to avail themselves of the new opportunity. In those times, this might have been partially due to the unreliability of data; the main problem was their general treatment. Then, statistical investigations are not reduced to mathematical calculations; circumstances accompanying the studied phenomena are also important, cf. Leibniz' opinion in § 4.1.1. Finally, their education did not prepare statisticians for grasping mathematical ideas and perhaps up to the 1870s they stubbornly held to *equally possible cases*, that is, to the theoretical probability.

Lack of such cases meant denial of the possibility to apply probability theory. But forget the 1870s! In 1916 A. A. Kaufman (Ploshko & Eliseeva 1990, p. 133) declared that the theory of probability is only applicable to independent trials with constant probability of *success* and certainly only when those equally possible cases existed.

Now, Quetelet. He had introduced mean inclinations to crime and marriage (although not for separate groups of population), but somehow statisticians did not for a long time understand that mean values ought not to be applied to individuals. As a consequence, by the end of his life and after his death (1874), mathematically ignorant statisticians went up in arms against those inclinations and probability in general (Rümelin 1867/1875, p. 25):

*Wenn mir die Statistik sagt, daß ich im Laufe des nächsten Jahres mit einer Wahrscheinlichkeit von 1 zu 49 sterben, mit einer noch*

*größeren Wahrscheinlichkeit schmerzliche Lücken in dem Kreis mir theurer Personen zu beklagen haben werde, so muß ich mich unter den Ernst dieser Wahrheit in Demuth beugen; wenn sie aber, auf ähnliche Durchschnittszahlen gestützt, mir sagen wollte, daß mit einer Wahrscheinlichkeit von 1 zu so und so viel* [I shall commit a crime] *so dürfte ich ihr unbedenklich antworten: ne sutor ultra crepidam*! [Cobbler! Stick to your last!].

A healthy man could have just as well rejected the conclusions drawn from a life table (Chuprov 1909/1959, pp. 211– 212).

Lexis infused a fresh spirit into (population) statistics. His followers, Bortkiewicz, Chuprov, Bohlmann, Markov, founded the so-called *continental direction of statistics*. In England, Galton, and Pearson somewhat later created the Biometric school which had been statistically studying Darwinism. The editors of its journal, *Biometrika*, a *Journal for the Statistical Study of Biological Problems*, were Weldon (a biologist who died in 1906), Pearson and Davenport (an author of a book on biometry and several articles) *in consultation* with Galton. Here is its Editorial published in 1902, in the first issue of that journal:

*The problem of evolution is a problem in statistics*. […] *We must turn to the mathematics of large numbers, to the theory of mass phenomena, to interpret safely our observations.* […] *May we not ask how it came about that the founder of our modern theory of descent made so little appeal to statistics?* […] *The characteristic bent of Darwin's mind led him to establish the theory of descent without mathematical conceptions; even so Faraday's mind worked in the case of electro-magnetism. But as every idea of Faraday allows of mathematical definition, and demands mathematical analysis*, […] *so every idea of Darwin – variation, natural selection* […] *– seems at once to fit itself to mathematical definition and to demand statistical analysis.* […] T*he biologist, the mathematician and the statistician have hitherto had widely differentiated field of work.* […] *The day will come* […] *when we shall find mathematicians who are competent biologists, and biologists who are competent mathematicians …*

During many years the Biometric school had been keeping to empiricism (Chuprov 1918 – 1919, t. 2, pp. 132 – 133) and he and Fisher (1922, pp. 311 and 329n) both indicated that Pearson confused theoretical and empirical indicators. And here is Kolmogorov (1947, p. 63; 1948/2002, p. 68):

*The modern period in the development of mathematical statistics began with the fundamental works of* […] *Pearson, Student, Fisher* […]. *Only in the contributions of the English school did the application of probability theory to statistics cease to be a collection of separate isolated problems and become a general theory of statistical testing of stochastic hypotheses* (*of hypotheses about laws of distribution*) *and of statistical estimation of parameters of these laws.*

*The main weakness of the* [Biometric] *school* [in 1912] *were: 1. Rigorous results on the proximity of empirical sample characteristics to the theoretical ones existed only for independent trials. 2. Notions of the logical structure of the theory of probability, which underlies all*

*the methods of mathematical statistics, remained at the level of eighteenth century results. 3. In spite of the immense work of compiling statistical tables […], in the intermediate cases between 'small' and 'large' samples their auxiliary techniques proved highly imperfect.*

I (2010) have collected many pronouncements about the Biometric school and Pearson; hardly known outside Russia was Bernstein's high opinion. I note that Kolmogorov passed over in silence the Continental direction of statistics. Chuprov had exerted serious efforts for bringing together that Continental direction and the Biometric school, but I am not sure that he had attained real success. And this I say in spite of the Resolution of condolence passed by the Royal Statistical Society after the death of its Honorary member (Sheynin 1990/2011, p. 156). It stated that Chuprov's contributions (not special efforts!) *did much to harmonize the methods of statistical research developed by continental and British workers*. Even much later Bauer (1955, p. 26) reported that he had investigated how both schools had been applying analysis of variance and concluded (p. 40) that their work was going on side by side but did not tend to unification.

## Bibliography
### Abbreviation

# Index of Names
The numbers refer to subsections rather than pages.